\RequirePackage{amsthm}

\documentclass[pdflatex,sn-mathphys-num]{sn-jnl}

\usepackage{graphicx}%
\usepackage{multirow}%
\usepackage{amsmath,amssymb,amsfonts}%
\usepackage{amsthm}%
\usepackage{mathrsfs}%
\usepackage[title]{appendix}%
\usepackage{xcolor}%
\usepackage{textcomp}%
\usepackage{manyfoot}%
\usepackage{booktabs}%
\usepackage{algorithm}%
\usepackage{algorithmicx}%
\usepackage{algpseudocode}%
\usepackage{listings}%

\usepackage{pgfgantt}
\usepackage{geometry}
\usepackage{relsize}
\usetikzlibrary{arrows.meta,arrows}
\usepackage{lmodern}
\usepackage{anyfontsize}
\usepackage[caption=false]{subfig}
\usepackage{adjustbox}
\usepackage{enumitem}
\usepackage{hhline}

\DeclareMathOperator{\EX}{\mathbb{E}}

\theoremstyle{thmstyleone}%
\theoremstyle{thmstyletwo}%

\theoremstyle{thmstylethree}%

\begin{document}

\title[Survey of a class of iterative row-action methods: The Kaczmarz method]{Survey of a class of iterative row-action methods: The Kaczmarz method}


\author*[1]{\fnm{Inês} \sur{A. Ferreira}}\email{ines.alves.ferreira@tecnico.ulisboa.pt} \equalcont{These authors contributed equally to this work.}

\author[1,2]{\fnm{Juan} \sur{A. Acebr\'on}}\email{juan.acebron@ist.utl.pt}
\equalcont{These authors contributed equally to this work.}

\author[1]{\fnm{Jos\'e} \sur{Monteiro}}\email{jcm@inesc-id.pt}
\equalcont{These authors contributed equally to this work.}

\affil*[1]{\orgdiv{INESC-ID}, \orgname{Instituto Superior Técnico}, \orgaddress{\city{Universidade de Lisboa}, \country{Portugal}}}

\affil[2]{\orgdiv{Department of Mathematics}, \orgname{Carlos III University of Madrid}, \orgaddress{\country{Spain}}}



\abstract{The Kaczmarz algorithm is an iterative method that solves linear systems of equations. It stands out among iterative algorithms when dealing with large systems for two reasons. First, at each iteration, the Kaczmarz algorithm uses a single equation, resulting in minimal computational work per iteration. Second, solving the entire system may only require the use of a small subset of the equations. These characteristics have attracted significant attention to the Kaczmarz algorithm. Researchers have observed that randomly choosing equations can improve the convergence rate of the algorithm. This insight led to the development of the Randomized Kaczmarz algorithm and, subsequently, several other variations emerged. In this paper, we extensively analyze the native Kaczmarz algorithm and many of its variations using large-scale dense random systems as benchmarks. Through our investigation, we have verified that, for consistent systems, various row sampling schemes can outperform both the original and Randomized Kaczmarz method. Specifically, sampling without replacement and using quasirandom numbers are the fastest techniques. However, for inconsistent systems, the Conjugate Gradient method for Least-Squares problems overcomes all variations of the Kaczmarz method for these types of systems.}




\maketitle

\section{Introduction}
\label{sec:intro}
Solving linear systems of equations is a fundamental problem in science and engineering. Therefore, it is important to create computer algorithms that can efficiently find solutions, especially for large-scale problems that generally require longer to solve.

Iterative methods that solve linear systems, unlike direct methods, do not produce an exact solution in a finite number of steps. Instead, they generate a sequence of approximate solutions that get increasingly closer to the exact solution with each iteration. When dealing with large systems, the computational cost of direct methods becomes excessive. For these systems, iterative methods are a good alternative that balances accuracy and computation time.

There are two special classes of iterative methods: row-action methods and column-action methods. These use a single row or column of the system's matrix per iteration. Each iteration of these methods has diminutive computational work compared with other iterative methods that use the entire matrix. An example of a row-action algorithm is the Kaczmarz algorithm, proposed by Kaczmarz himself in 1937 \cite{kaczmarz1937angenaherte} and the focus of this paper. Methods like Kaczmarz that use one equation at a time can be used in real-time while data is being collected. Furthermore, if the amount of data is so large that storing it in a machine is not possible, row/column action methods are still able to solve the system.

An additional issue is that some problems may have multiple solutions, and others none. One example of the application of solving linear systems in the real-world are problems derived from computed tomography (CT). During a CT scan, one has to reconstruct an image of the scanned body using radiation data measured by detectors. Row/column action methods are used to solve problems derived from CT scans to the detriment of other iterative algorithms since they can process information while data is being collected. Additionally, raw data obtained from a CT scan demands substantial memory resources and might surpass the storage capacity of a single machine, especially considering that some accurate 3D CT scans may generate hundreds of gigabytes of data per second \cite{bicer2017trace}. Consequently, methods that require processing the entire dataset at once cannot be used. Moreover, physical quantities are always measured with some error and, since CT data is not an exception, the problem of reconstructing images is hampered by noise, leading to a system with no solution.

In this paper, we analyze and compare the performance of several Kaczmarz-based methods that sample matrix rows according to different criteria. More specifically, we evaluate the relationship between the number of iterations needed to find a solution and the corresponding execution time, for systems with different sizes.

The organization of this document is as follows. In Section~\ref{sec:theory}, we describe several types of linear systems and respective solutions. In Section~\ref{sec:seq_var} we introduce the original algorithm together with some of the methods inspired by the Randomized Kaczmarz algorithm. In Section~\ref{sec:applic_kacz} we discuss some applications of the Kaczmarz method including adaptations for linear systems of inequalities and linear systems from CT scans. In Section~\ref{sec:res_seq}, we present the implementation details of the algorithm and some of its variations, together with experimental results. Finally, in Section~\ref{sec:conclusion}, we conclude this paper and discuss future work.

\section{Background}
\label{sec:theory}
We start by outlining the categorization of linear systems based on factors such as matrix dimensions and the existence of a solution, followed by a small introduction to solving linear systems. Then we introduce the Kaczmarz algorithm and its key attributes. In the remainder of this section, we present several iterative methods, some of which are modifications of the Kaczmarz method.

\subsection{Types of Linear Systems}

A linear system of equations can be written as
\begin{equation} \label{eq:system}
    Ax = b \: ,
\end{equation}
where $A \in \mathbb{R}^{m \times n}$ is a matrix of real elements, $x \in \mathbb{R}^{n}$ is called the solution and $b \in \mathbb{R}^{m}$ is a vector of constants. Linear systems can be classified based on distinct criteria.

When we consider the existence of a solution, systems can be categorized as either \textbf{consistent} or \textbf{inconsistent}. Let $[ A\,|\,b ]$ be the augmented matrix that represents the system (\ref{eq:system}). Remember that the rank of a matrix corresponds to the number of linearly independent rows/columns and that a full rank matrix is a matrix that has the maximum possible rank, given by $min(m,n)$. A system is \textbf{consistent} if there is at least one value of $x$ that satisfies all the equations in the system. To check for consistency one can compute the rank of the original and augmented matrix. If $rank[A]$ is equal to $rank[A|b]$, then the system is consistent. For consistent systems, if $rank[A]=rank[A|b]=min(n,m)$, there is a single solution that satisfies (\ref{eq:system}). Otherwise, if $rank[A]=rank[A|b]<min(n,m)$, matrix $A$ is rank-deficient and there are infinite solutions. However, it is possible that there are equations that contradict each other and, in that case, the system is said \textbf{inconsistent}. In this case, there is no solution, and $rank[A]<rank[A|b]$.

When considering the relationship between the number of equations and variables, systems can be classified as \textbf{overdetermined} or \textbf{underdetermined}. When $m \geq n$, there are more equations than variables, and the system is said to be \textbf{overdetermined}. Note that, for these systems, a full rank matrix has rank $min(m,n)=n$. If the system is consistent and there is a single solution that satisfies (\ref{eq:system}), we denote that solution by $x^*$. If the system is inconsistent, we are usually interested in finding the least-squares solution, that is
\begin{equation} \label{eq:ls_system}
    x_{LS} = arg \min_{x} \| Ax - b \|^2 \:.
\end{equation}
Throughout this document, we use $\| \: . \: \|$ to represent the Euclidean $L^2$ norm. In real-world overdetermined systems, it is more frequent to have inconsistent systems than consistent systems. The least-squares solution can be computed using the normal equations $x_{LS} = (A^{T} A)^{-1} A^{T} b \: = A^{\dagger} b$, where $A^{\dagger}$ is the Moore–Penrose inverse or pseudoinverse. If the system has a matrix with $m = n$, we are in the presence of an exactly determined system.

\textbf{Underdetermined} systems verify $m < n$, meaning that we have fewer equations than variables. When an underdetermined system is consistent, even if it is full rank, there will still be $n - m$ degrees of freedom and the system has infinite solutions. In this case, we are often interested in the least Euclidean norm solution,
\begin{equation} \label{eq:sys_under}
x_{LN} = arg \min_{x} \: \|x\| \quad \text{subject to} \quad Ax = b \: ,
\end{equation}
which can be computed using $x_{LN} = A^{T} (A A^{T})^{-1} b$.

\subsection{Solving Linear Systems}

There are two classes of numerical methods that solve linear systems of equations: direct and indirect (or iterative). Direct methods compute the solution of the system in a finite number of steps. Although solutions given by direct methods have a high level of precision, the computational cost and memory usage can be high for large matrices when using these methods. Precision can be defined as the Euclidean distance of the difference between the obtained solution and the real solution of the system. Some widely used direct methods are Gaussian Elimination and LU Factorization. Iterative methods generate a sequence of approximate solutions that get increasingly closer to the exact solution with each iteration. Iterative methods are generally less computationally demanding than direct methods and their convergence rate is highly dependent on the properties of the matrix of the system. The rate of convergence \cite{senning2007computing} of an algorithm quantifies how quickly a sequence approaches its limit. In the case of iterative methods, a higher rate of convergence means that the method will take fewer iterations to approach the solution. The precision of the solutions given by iterative methods can be controlled by external parameters, meaning that there is a trade-off between the running time and the desired accuracy of the solution. Therefore, if high precision is not a requirement, iterative methods can outperform direct methods, especially for large and sparse matrices. Some popular iterative methods are the Jacobi method and the Conjugate Gradient (CG) method.

A particular class of iterative methods is row-action methods. These use only one row of matrix $A$ in each iteration, meaning that the computational work per iteration is small compared to methods that make use of the entire matrix. There are also column-action methods that, instead of using a single row per iteration, use a single column. An example of a row action algorithm is the Kaczmarz method which is introduced next.


\section{Variations of the Kaczmarz Method} \label{sec:seq_var}

In this section, we introduce the Kaczmarz method and several of its variations. We also introduce some other randomized iterative methods that were inspired by the development of the Randomized Kaczmarz method. We finish the section with a small introduction to some parallelization strategies for the original and randomized Kaczmarz methods.

\subsection{The Kaczmarz Method} \label{sec:CK}

The Kaczmarz method \cite{kaczmarz1937angenaherte} is an iterative algorithm that solves consistent linear systems of equations $A x = b$. Let $A^{(i)}$ be the $i$-th row of $A$ and $b_i$ be the $i$-th coordinate of $b$. The original version of the algorithm can then be written as
\begin{equation} \label{eq:alg}
    x^{(k+1)} = x^{(k)} + \alpha_i \: \frac{b_i \: - \langle A^{(i)}, x^{(k)}  \rangle}{\|A^{(i)}\|^2} \: {A^{(i)}}^T \: , \quad \text{with} \quad i = k \text{ mod } m \: ,
\end{equation}
with $k$, the iteration number, starting at $0$ and where $\alpha_i$ is a relaxation parameter that can take any value between 0 and 2. Throughout this document we will use $\alpha_i = 1$ since most variations of the Kaczmarz method discussed here use this convention. Furthermore, the initial guess $x^{(0)}$ is typically set to zero, and $\langle , \rangle$ represents the dot product of two vectors. At each iteration, the estimated solution will satisfy a different constraint, until a point is reached where all the constraints are satisfied. Since the rows of matrix $A$ are used cyclically, this algorithm is also known as the Cyclic Kaczmarz (CK) method. The Kaczmarz method also converges to $x_{LN}$ in the case of underdetermined systems. For more details see Section 3.3 of \cite{ma2015convergence}.

\begin{figure}[t]
    \centering
    \begin{minipage}{.5\textwidth}\centering
    \subfloat[Example of a consistent system with an unique solution.]{\includegraphics[width=0.95\columnwidth]{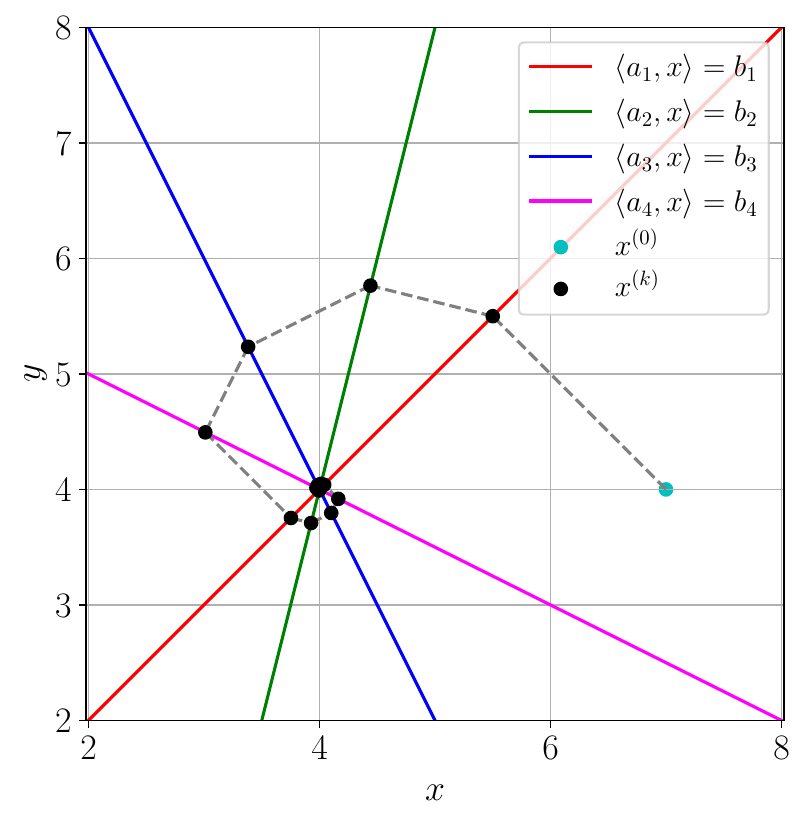}\label{fig:example_cyclic_consist}
    }
    \end{minipage}%
    \begin{minipage}{.5\textwidth}\centering
    \subfloat[Example of an inconsistent system and its least-squares solution.]{\includegraphics[width=0.95\columnwidth]{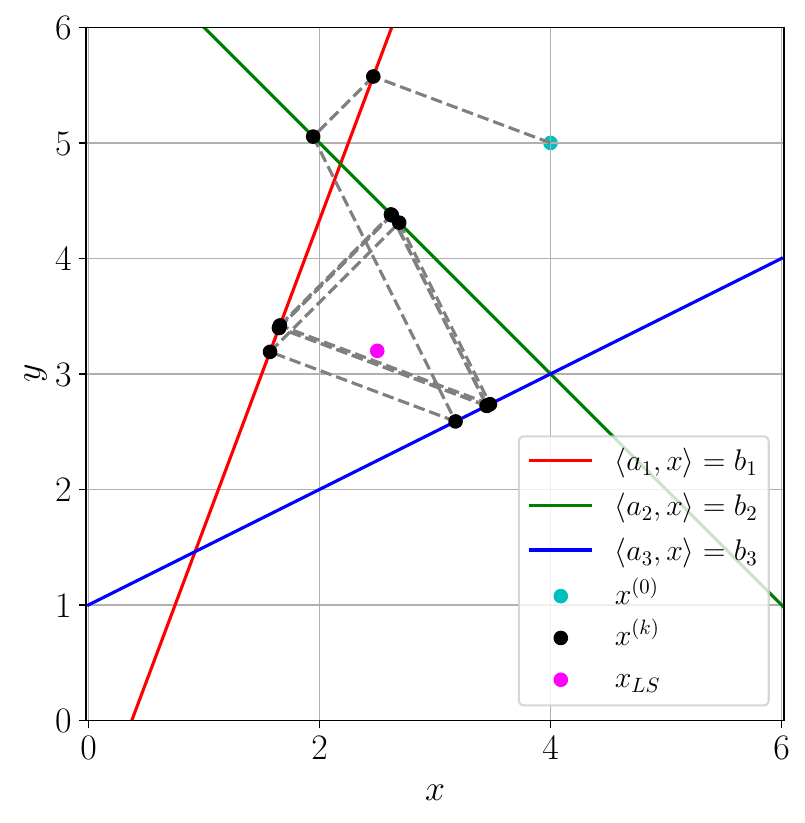}\label{fig:example_cyclic_inconsist}
    }
    \end{minipage}
    \caption{Convergence of the Kaczmarz method in 2 dimensions for consistent and inconsistent systems.}
    \label{fig:example_cyclic}
\end{figure}

The Kaczmarz method has a clear geometric interpretation (if the relaxation parameter $\alpha_i$ is set to 1). Each iteration $x^{(k+1)}$ can be thought of as the projection of $x^{(k)}$ onto $H_i$, the hyperplane defined by $H_i = \{x : \langle A^{(i)}, x \rangle = b_i \}$. Figure~\ref{fig:example_cyclic} provides the geometrical interpretation of a linear system with 4 equations in 2 dimensions: each equation of the system is represented by a line in a plane; the evolution of the estimate of the solution, $x^{(k)}$, is also shown throughout several iterations. Since we are using the rows of the matrix cyclically, we start by projecting the initial estimate of the solution onto the first equation and then the second, and so on. Figure~\ref{fig:example_cyclic_consist} shows a consistent system with a unique solution, $x^*$, the point where all equations meet. It is clear that the estimate of the solution given by the Kaczmarz method will eventually converge to some vector $x^*$. In Figure~\ref{fig:example_cyclic_inconsist} we have an inconsistent system with its least-squares solution marked by the magenta dot. Note that Kaczmarz can't reach it since the estimate of the solution will always remain a certain distance from $x_{LS}$. In the next section, we will show how this distance can be quantified.

\begin{figure}[t]
    \centering
    \begin{minipage}{.5\textwidth}\centering
    \subfloat[Cyclical selection of rows.]{\includegraphics[width=0.95\columnwidth]{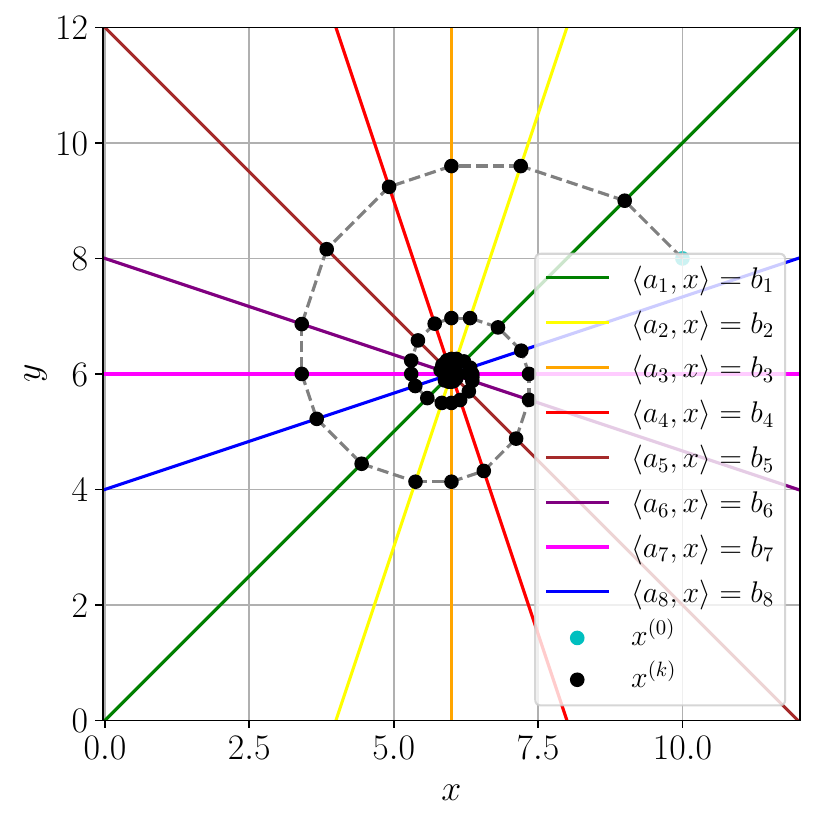}\label{fig:example_cyclic_coherent}}
    \end{minipage}%
    \begin{minipage}{.5\textwidth}\centering
    \subfloat[Random selection of rows.]{\includegraphics[width=0.95\columnwidth]{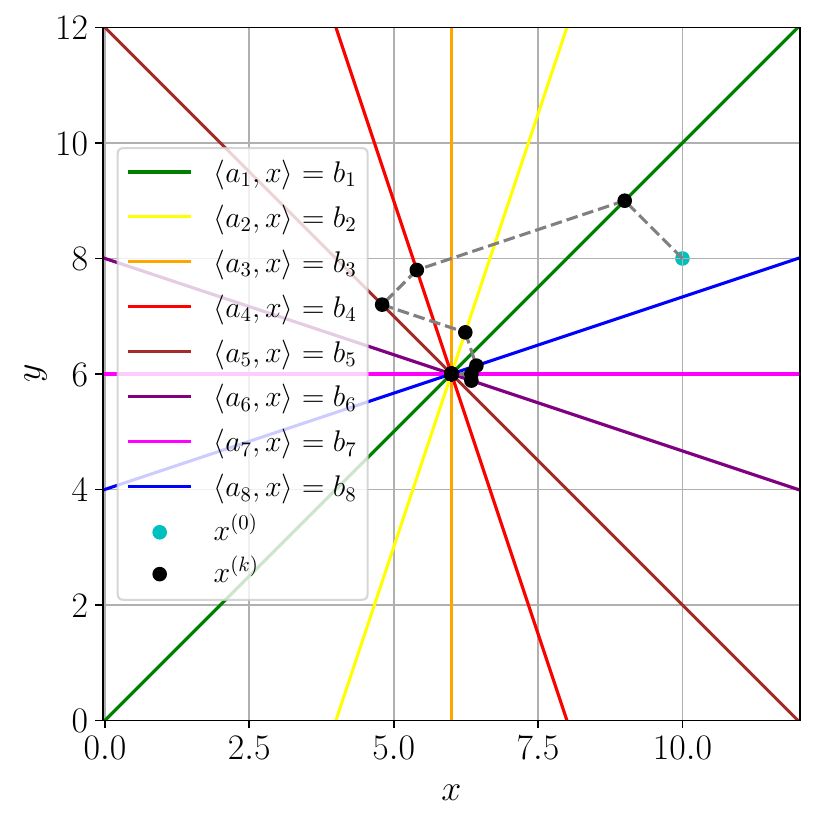}\label{fig:example_rand_coherent}
    }
    \end{minipage}
    \caption{Convergence of the Kaczmarz method for a consistent system in 2 dimensions using two different row selection criteria.}
    \label{fig:example_coherent}
\end{figure}

Figure~\ref{fig:example_coherent} shows the special case of a consistent system with a highly coherent matrix, defined as a matrix for which the angle between consecutive rows is small. If we go through the matrix cyclically, as shown in Figure~\ref{fig:example_cyclic_coherent}, the convergence is quite slow since the change in the solution estimate is minimal. On the other hand, in Figure~\ref{fig:example_rand_coherent}, where the rows of the matrix are used in a random fashion, the estimate of the solution approaches the system solution much faster. This was experimentally observed and led to the development of randomized versions of the Kaczmarz method.

\subsection{Randomized Kaczmarz Method} \label{sec:RK}

Kaczmarz \cite{kaczmarz1937angenaherte} showed that, if the linear system is consistent with a unique solution, the method converges to $x^*$, but the rate at which the method converges is very difficult to quantify. Known estimates for the rate of convergence of the cyclic Kaczmarz method rely on quantities of matrix $A$ that are hard to compute and this makes it difficult to compare with other algorithms.
Ideally, we would like to have a rate of convergence that depends on the condition number of the matrix $A$, $k(A) = \|A\| \: \|A^{-1}\|$. Since $\|A\| = \sigma_{max}(A)$ and $\|A^{-1}\| = 1/\sigma_{min}(A)$, the condition number of $A$ can also be written as $k(A) = \sigma_{max}(A)/\sigma_{min}(A)$. Finding a rate of convergence that depends on $k(A)$ is difficult for this algorithm since it relies on the order of the rows of the matrix, and the condition number is only related to the geometric properties of the matrix.

It has been observed 
that, instead of using the rows of $A$ in a cyclic manner, choosing rows randomly can improve the rate of convergence of the algorithm. To tackle the problem of finding an adequate rate of convergence for the Kaczmarz method and to try to explain the empirical evidence that randomization accelerates convergence, Strohmer and Vershynin \cite{Strohmer2007ARK} introduced a randomized version of the Kaczmarz method that converges to $x^*$. Instead of selecting rows cyclically, in the randomized version, in each iteration, we use the row with index $i$, chosen at random from the probability distribution
\begin{equation} \label{eq:prob_line}
    P\{ i = l \} = \frac{\|A^{(l)}\|^2}{\|A\|_F^2} \quad (l = 0, 1, 2, ..., m-1) \: ,
\end{equation}
where the Frobenius norm of a matrix is given by $\|A\|_F = \sqrt{\sum_{i=1}^{m} \|A^{(i)}\|^{2}}$. This version of the algorithm is called the Randomized Kaczmarz (RK) method. Strohmer and Vershynin proved that RK has exponential error decay, also known as linear convergence~\footnote{The concept of linear convergence in numerical analysis is often referred to, by mathematicians, as exponential convergence.}. Let $x_0$ be the initial guess, $x^*$ be the solution of the system, $\sigma_{min}(A)$ be the smallest singular value of matrix $A$, and $\kappa(A) = \|A\|_F \: \|A^{-1}\| = \|A\|_F / \sigma_{min}(A)$ be the scaled condition number of the system's matrix. The expected convergence of the algorithm for a consistent system can then be written as
\begin{equation} \label{eq:conv}
    \EX \| x^* - x^{(k)} \|^2 \leq ( 1 - \kappa(A)^{-2} )^k \| x^* - x^{(0)} \|^2 \leq \Big ( 1 - \frac{\sigma_{min}^2(A)}{\|A\|_F^2} \Big )^k \| x^* - x^{(0)} \|^2 \: .
\end{equation}
Later, Needell \cite{needell2010randomized} extended these results for inconsistent systems, showing that RK reaches an estimate that is within a fixed distance from the solution. This error threshold is dependent on the matrix $A$ and can be reached with the same rate as in the error-free case. The expected convergence proved by Strohmer and Vershynin can then be extended such that  
\begin{equation}
    \EX \| x^* - x^{(k)} \|^2 \leq \Big ( 1 - \frac{\sigma_{min}^2(A)}{\|A\|_F^2} \Big )^k \| x^* - x^{(0)} \|^2 \: + \frac{\|r_{LS}\|^2}{\sigma_{min}^2(A)} \: ,
\end{equation}
where $r_{LS} = b - A x_{LS}$ is the least-squares residual. This extra term is called the convergence horizon and is zero when the linear system is consistent.
Note that the rate of convergence of the algorithm depends only on the scaled condition number of $A$, and not on the number of equations in the system. This analysis shows that it is not necessary to know the whole system to solve it, but only a small part of it. Strohmer and Vershynin show that in extremely overdetermined systems the Randomized Kaczmarz method outperforms all other known algorithms and, for moderately overdetermined systems, it outperforms the celebrated conjugate gradient method. They reignited not only the research on the Kaczmarz method but also triggered the investigation into developing randomized linear solvers.

\subsection{Simple Randomized Kaczmarz Method} \label{sec:SRK}

Apart from the comparison between the original Kaczmarz method with the Randomized Kaczmarz method, Strohmer and Vershynin \cite{Strohmer2007ARK} also compared both these methods with the Simple Randomized Kaczmarz (SRK) method. In this method, instead of sampling rows using their norms, rows are sampled using a uniform probability distribution. Although they did not prove the convergence rate for this method, Schmidt~\cite{notes_on} did:
\begin{equation}
    \EX \| x^* - x^{(k)} \|^2 \leq \Big ( 1 - \frac{\sigma_{min}^2(A)}{m\|A\|_{\infty}^2} \Big )^k \| x^* - x^{(0)} \|^2.
\end{equation}
Schmidt also mentions that RK should be at least as fast as SRK, and faster if any two rows don’t have the same norm.

The work by Strohmer and Vershynin in \cite{Strohmer2007ARK} motivated other developments in row/column action methods by randomizing classical algorithms. In the following sections, we present several iterative methods, some of which are modifications to the RK method.

\subsection{Sampling Rows Using Quasirandom Numbers} \label{sec:quasirand_intro}

\begin{figure}[t]
    \centering
    \subfloat[Random numbers generated using an uniform probability distribution.]{\includegraphics[width=0.6\columnwidth]{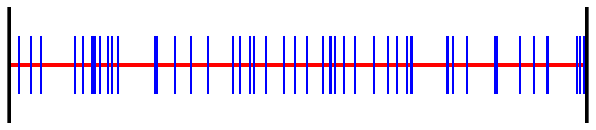}\label{fig:quasi_rand_dist_1}} \\
    \subfloat[Quasirandom numbers generated using the Sobol sequence.]{\includegraphics[width=0.6\columnwidth]{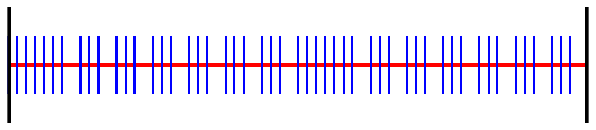}\label{fig:quasi_rand_dist_2}} \\
    \subfloat[Quasirandom numbers generated using the Halton sequence.]{\includegraphics[width=0.6\columnwidth]{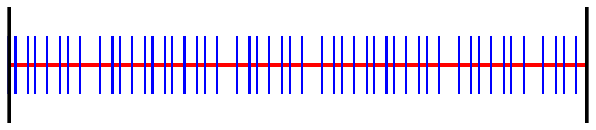}\label{fig:quasi_rand_dist_3}}
    \vspace{10pt}
    \caption{Distribution of 50 sampled numbers in interval $\left[ 1, 1000 \right] $.}
    \label{fig:quasi_rand_dist}
\end{figure}

In Section~\ref{sec:CK} we discussed how sampling rows in a random fashion instead of cyclically can improve convergence for highly coherent matrices. Furthermore, convergence should be faster if indices corresponding to consecutively sampled rows are not close so that the angle between these rows is not small. However, when sampling random numbers, these can form clumps and they can be poorly distributed, meaning that even if we sample rows randomly, these can still be very close regarding their position in the matrix and, consecutively, in the angle between them. To illustrate the importance of how the generation of random numbers affects the results of the Kaczmarz algorithm, Figure~\ref{fig:quasi_rand_dist_1} shows 50 random numbers sampled from the interval $\left[ 1, 1000 \right]$ using a uniform probability distribution, illustrating the process of sampling row indices for a matrix with 1000 rows. Note that there are, simultaneously, clumps of numbers and areas with no sampled numbers.

This is where quasirandom numbers, also known as low-discrepancy sequences, come in: they are sequences of numbers that are evenly distributed, meaning that there are no areas with a very low or very high density of sampled numbers. Quasirandom numbers have been shown to improve the convergence rate of Monte-Carlo-based methods, namely methods related to numerical integration \cite{sloan1994lattice}.

Several low-discrepancy sequences can be used to generate quasirandom numbers. Here we will work with two sequences that are widely used: the Halton sequence~\cite{halton1960efficiency} and the Sobol sequence \cite{sobol1976uniformly}. To show that indeed these sequences can generate numbers that are more evenly distributed than the ones generated using pseudo-random numbers, Figures~\ref{fig:quasi_rand_dist_2} and~\ref{fig:quasi_rand_dist_3} show 50 sampled numbers from the interval $\left[ 1, 1000 \right] $ using the Sobol and Halton sequences. We will denote the two methods of using quasirandom numbers for row sampling as SRK-Halton and SRK-Sobol.



\subsection{Randomized Block Kaczmarz Method} \label{sec:RBK}

In the Kaczmarz method, a single constraint is enforced in each iteration. To speed up the convergence of this method, a block version \cite{elfving1980block} was designed such that several constraints are enforced at once. Needell and Tropp \cite{needell2014paved} were the first to develop a randomized block version of the Kaczmarz method, the Randomized Block Kaczmarz (RBK) method, and to provide its rate of convergence. Each block is a subset of rows of matrix $A$ chosen with a randomized scheme. Let us consider that, in each iteration $k$, we select a subset of row indices of $A$, $\tau_k$. The current iterate $x^{(k+1)}$ is computed by projecting the previous iterate $x^{(k)}$ onto the solution of $A_{\tau_k} x = b_{\tau_k}$, that is
\begin{equation} \label{eq:rbk}
    x^{(k+1)} = x^{(k)}+(A_{\tau_k})^{\dagger} (b_{\tau_k} - A_{\tau_k} x^{(k)}) \: ,
\end{equation}
where $(A_{\tau_k})^{\dagger}$ denotes the pseudoinverse. The subset of rows $\tau_k$ is chosen using the following procedure. First, we define the number of blocks, denoted by $l$. Then we divide the rows of the matrix into $l$ subsets, creating a partition $T = \{ \tau_1, ..., \tau_l \}$. The block $\tau_k$ can be chosen from the partition using one of two methods: it can be chosen randomly from the partition independently of all previous choices, or it can be sampled without replacement, an alternative that Needell and Tropp found to be more effective. In the latter, the block can only contain rows that have yet to be selected. Only when all rows are selected can there be a re-usage of rows. The partition of the matrix into different blocks is called row-paving. If $l = m$, that is, if the number of blocks is equal to the number of rows, each block consists of a single row and we recover the RK method. The calculation of the pseudoinverse $(A_{\tau})^{\dagger}$ in each iteration is a computationally expensive step. But, if the submatrix $A_{\tau}$ is well-conditioned, we can use algorithms like Conjugate Gradient for Least-Squares (CGLS) to efficiently calculate it.

To analyze the rate of convergence of this version some quantities need to be defined since they characterize the performance of the algorithm. The row paving $(l, \alpha, \beta)$ of a matrix $A$ is a partition $T = \{ \tau_1, ..., \tau_l \}$ that verifies
\begin{equation}
    \alpha \leq \lambda_{min} (A_{\tau} A_{\tau}^*) \quad \text{and} \quad \lambda_{max} (A_{\tau} A_{\tau}^*) \leq \beta \quad \text{for each} \quad \tau \in T \: ,
\end{equation}
where $l$ is the number of blocks and describes the size of the paving and the quantities $\alpha$ and $\beta$ are called lower and upper paving bounds. Suppose $A$ is a matrix with full column rank and row paving $(l, \alpha, \beta)$. Then, the rate of convergence of RBK is described by
\begin{equation}
    \EX \|x^{(k)}-x^*\|^2 \leq \Big ( 1 - \frac{\sigma_{min}^2(A)}{\beta l}\Big )^{k} \|x^{(0)}-x^*\|^2 + \frac{\beta}{\alpha} \frac{\|r\|^2}{\sigma_{min}^2(A)} \: ,
\end{equation}
where $r = Ax - b$ is the residual.
Similarly to RK, the Randomized Block Kaczmarz method exhibits an expected linear rate of convergence. Notice that the rate of convergence for a consistent system only depends on the upper paving bound $\beta$ and on the paving size $l$. Regarding inconsistent systems, the convergence horizon is also dependent on $\beta / \alpha$, denoted the paving conditioning. In summary, the paving of a matrix is intimately connected with the rate of convergence.

It is important to determine in which cases it is beneficial to use the RBK method to the detriment of the Kaczmarz method. The first case is if we have matrices with good row paving, that is, the blocks are well conditioned. In the presence of good row paving, the matrix-vector multiplication and the computation of $(A_{\tau_k})^{\dagger}$ can be very efficient, meaning that one iteration of the block method has roughly the same cost as one iteration of the standard method. Depending on the characteristics of matrix $A$, there can be efficient computational methods to produce good paving if the matrix does not naturally have good paving. A second case is connected with the implementation of the algorithm: the simple Kaczmarz method transfers a new equation into working memory in each iteration, making the total time spent in data transfer throughout the algorithm excessive; on the other hand, not only does the block Kaczmarz algorithm move large data blocks into working memory at a time, but it also relies on matrix-vector multiplication that can be accelerated using basic linear algebra subroutines (BLASx) \cite{wang2016blasx}.

\subsection{Randomized Coordinate Descent Method} \label{sec:RGS}

Leventhal and Lewis \cite{leventhal2010randomized} developed the Randomized Coordinate Descent Method, also known as the Randomized Gauss-Seidel (RGS) algorithm. Like RK, for overdetermined consistent systems, this algorithm converges to the unique solution $x^*$. Unlike RK, for overdetermined inconsistent systems, RGS converges to the least squares solution $x_{LS}$; for underdetermined consistent systems RGS does not converge to the least euclidean norm solution $x_{LN}$ \cite{ma2015convergence}.
RGS is an iterative method that uses only one column of matrix $A$ in each iteration, represented by $A_{(l)}$,  which is chosen at random from the probability distribution
\begin{equation} \label{eq:prob_cols}
    P\{ j = l \} = \frac{\|A_{(l)}\|^2}{\|A\|_F^2} \quad (l = 0, 1, 2, ..., n-1) \: .
\end{equation}
It also uses an intermediate variable $r \in \mathbb{R}^{m}$. The algorithm has three steps to be completed in a single iteration:
\begin{equation}
    \alpha^{(k)} = \frac{\langle A_{(j)}, r^{(k)} \rangle}{\|A_{(j)}\|^2} \: ,
\end{equation}
\begin{equation}
    x^{(k+1)} = x^{(k)} + \alpha^{(k)} e_{(j)} \: , \quad \text{and}
\end{equation}
\begin{equation}
    r^{(k+1)} = r^{(k)} - \alpha^{(k)} A_{(j)} \: ,
\end{equation}
where $e_{(j)}$ is the jth coordinate basis column vector (all zeros with a 1 in the $j$th position) and $r$ is initialized to $b - A x^{(0)}$.
Leventhal and Lewis also showed that, like RK, RGS converges linearly in expectation.
Later, Needell, Zhao, and Zouzias~\cite{needell2015randomized} developed a block version of the RGS method, called the Randomized Block Coordinate Descent (RBGS) method, that uses several columns in each iteration, with the goal of improving the convergence rate. 

\subsection{Randomized Extended Kaczmarz Method} \label{sec:REK}

The Randomized Kaczmarz method can only be applied to consistent linear systems, but most systems in real-world applications are affected by noise, and are consequently inconsistent, meaning that it is important to develop a version of RK that solves least-squares problems.

To extend the work developed by Strohmer and Vershynin \cite{Strohmer2007ARK} to inconsistent systems, Zouzias and Freris \cite{zouzias2013randomized} introduced the Randomized Extended Kaczmarz (REK) method.
This method is a combination of the Randomized Orthogonal Projection algorithm \cite{rokhlin2009fast} together with the Randomized Kaczmarz method and it converges linearly in expectation to the least-squares solution, $x_{LS}$. The algorithm is a mixture of a row and column action method since, in each iteration, we use one row and one column of matrix $A$. Rows are chosen with probability proportional to the row norms and columns are chosen with probability proportional to column norms. In each iteration, we use a row with index $i$, chosen from the probability distribution in (\ref{eq:prob_line}) and a column with index $j$, chosen at random from the probability distribution (\ref{eq:prob_cols}). Each iteration of the algorithm can then be computed in two steps. In the first one, the projection step, we calculate the auxiliary variable $z \in \mathbb{R}^{m}$. In the second one, the Kaczmarz step, we use $z$ to estimate the solution $x$. Variables $z$ and $x$ are initialized, respectively, to $b$ and $0$. The two steps of each iteration are then
\begin{equation}
    z^{(k+1)} = z^{(k)} - \frac{\langle A_{(j)}, z^{(k)} \rangle}{\|A_{(j)}\|^2} \: A_{(j)} \: , \quad \text{and}
\end{equation}
\begin{equation}
    x^{(k+1)} = x^{(k)} + \frac{b_i - z_i^{(k)} - \langle x^{(k)}, A^{(i)} \rangle}{\|A^{(i)}\|^2} \: {A^{(i)}}^T \: .
\end{equation}

The rate of convergence can be described by
\begin{equation}
    \EX \| x_{LS} - x^{(k)} \|^2 \leq \Big( 1 - \frac{\sigma_{min}^2(A)}{\|A\|_F^2} \Big)^{ \lfloor k/2 \rfloor} \: \Big(  1 + 2 \:  \frac{\sigma_{max}^2(A)}{\sigma_{min}^2(A)} \Big) \: \| x_{LS} \|^2 \: .
\end{equation}


\subsection{Randomized Double Block Kaczmarz Method} \label{sec:RDBK}

We have seen in the previous section that the REK method is a variation of the RK method that converges to the least squares solution $x_{LS}$. In Section~\ref{sec:RBK} we described situations where the RBK method for consistent systems can outperform the non blocked version. Needell, Zhao, and Zouzias \cite{needell2015randomized} developed a new algorithm called the Randomized Double Block Kaczmarz (RDBK) method that combines both the REK and the RBK algorithm to develop a method that solves inconsistent systems with accelerated convergence. Recall that in REK there is a projection step that makes use of a single column of matrix $A$ and there is the Kaczmarz step that utilizes a single row. This means that, in a block version of REK, we need both a column partition for the projection step and a row partition for the Kaczmarz step, hence the name ''Double Block´´. In each iteration $k$ we select a subset of row indices of $A$, $\tau_k$, and a subset of column indices of $A$, $v_k$. The projection step and the Kaczmarz step can then be written as
\begin{equation}
        z^{(k+1)} = z^{(k)} - A_{v_k} (A_{v_k})^{\dagger} z^{(k)} \: , \quad \text{and}
\end{equation}
\begin{equation}
        x^{(k+1)} = x^{(k)} + (A_{\tau_k})^{\dagger} (b_{\tau_k} - z_{\tau_k}^{(k+1)} - A_{\tau_k} x^{(k)}) \: .
\end{equation}
The initialization variables are $x^{(0)} = 0$ and $z^{(0)} = b$.
Let $l$ be the number of row blocks, $\overline{l}$ be the number of column blocks, $\alpha$ and $\beta$ be the lower and upper paving row bounds and $\overline{\alpha}$ and $\overline{b}$ be the lower and upper paving column bounds. The row paving will be described by $(l, \alpha, \beta)$ and the column paving by $(\overline{l}, \overline{\alpha}, \overline{\beta})$ and the rate of convergence is represented by
\begin{equation}
    \EX \|x^{(k)}-x_{LS}\|^2 \leq \gamma^k \|x^{(0)}-x_{LS}\|^2 + \Big ( \gamma^{\lfloor k/2 \rfloor} + \overline{\gamma}^{\lfloor k/2 \rfloor} \Big ) \frac{\| b_{\mathcal{R}(A)} \|^2}{\alpha (1 - \gamma)} \: ,
\end{equation}
where $b_{\mathcal{R}(A)}$ denotes the projection of $b$ onto the row space of matrix $A$, $\gamma = 1 - \sigma_{min}^2(A)/(l \beta)$ and $\overline{\gamma} = 1 - \sigma_{min}^2(A)/(\overline{l} \overline{\beta})$.

Just like for RBK, the use of the RDBK to the detriment of the REK method depends on the row/column paving of the matrix.  In addition to the block version of REK, the authors also show that the RGS method can be accelerated using a block version (RBGS) (see Section~\ref{sec:RGS}). However, RBGS has an advantage regarding the RBK method: the RBGS method only requires column paving since the RGS method uses a single column per iteration.


\subsection{Greedy Randomized Kaczmarz Method} \label{sec:grk}

The Greedy Randomized Kaczmarz (GRK) method introduced by Bai and Wu \cite{bai2018greedy} is a variation of the Randomized Kaczmarz method with a different row selection criterion. Note that the selection criterion for rows in the RK method can be simplified to uniform sampling if we scale matrix $A$ with a diagonal matrix that normalizes the Euclidean norms of all its rows. But, in iteration $k$, if the residual vector $r^{(k)} = b - A x^{(k)}$ has $|r^{(k)}_i| > |r^{(k)}_j|$, we would like for row $i$ to be selected with a higher probability than row $j$. In summary, GRK differs from RK by selecting rows with larger entries of the residual vector with higher probability. Each iteration is still calculated using (\ref{eq:alg}) but the row index $i$ chosen in iteration $k$ is computed using the following steps:

Compute
\begin{equation}
    \epsilon_k = \frac{1}{2} \Big( \frac{1}{\|b - A x^{(k)}\|^2} \max_{1 \leq i \leq m} \Big \{ \frac{|b_i - \langle x^{(k)}, A^{(i)} \rangle|^2}{\| A^{(i)} \|^2} \Big \} + \frac{1}{\|A\|_F^2} \Big ) \: .
\end{equation}
Determine the index set of positive integers
\begin{equation}
    \quad \mathcal{U}_k = \Big \{ i : | b_i - \langle x^{(k)}, A^{(i)} \rangle |^2 \geq \epsilon_k \|b - A x^{(k)}\|^2 \: \| A^{(i)} \|^2 \Big \} \: .
\end{equation}
Compute vector
\begin{equation}
    \tilde{r}^{(k)}_{i} =
    \begin{cases}
      b_i - \langle x^{(k)}, A^{(i)} \rangle, & \text{if $i \in \mathcal{U}_k$} \\
      0, & \text{otherwise.}
    \end{cases}
\end{equation}
Select $i_k \in \mathcal{U}_k$ with probability
\begin{equation}
    P\{i_k = i\} = \frac{| \tilde{r}^{(k)}_{i} |^2}{\| \tilde{r}^{(k)} \|^2} \: .
\end{equation}
The Greedy Randomized Kaczmarz method presents a faster convergence rate when compared to the Randomized Kaczmarz method, meaning that it is expected for GRK to outperform RK.

\subsection{Selectable Set Randomized Kaczmarz Method} \label{ssrk}

Just as the Greedy Randomized Kaczmarz method, the Selectable Set Randomized Kaczmarz (SSRK) method \cite{yaniv2021selectable} is a variation of the Randomized Kaczmarz method with a different probability criterion for row selection that avoids sampling equations that are already solved by the current iterate. The set of equations that are not yet solved is referred to as the selectable set, $\mathcal{S}_k$. In each iteration of the algorithm, the row to be used in the Kaczmarz step is chosen from the selectable set. The algorithm is as follows
\begin{itemize}
    \item Sample row $i_k$ according to probability distribution in (\ref{eq:prob_line}) with rejection until $i_k \in \mathcal{S}_k$;
    \item Compute the estimate of the solution $x^{(k+1)}$;
    \item Update the selectable set $\mathcal{S}_{k+1}$ such that $i_k \notin \mathcal{S}_{k+1}$.
\end{itemize}
The last step of the algorithm can be done in two ways, meaning that there are two variations of the SSRK method: the Non-Repetitive Selectable Set Randomized Kaczmarz (NSSRK) method and the Gramian Selectable Set Randomized Kaczmarz (GSSRK) method. 

For the NSSRK method, the selectable set is updated in each iteration by simply using all rows except the row that was sampled in the previous iteration. In that case, the last step of the algorithm can be written as $\mathcal{S}_{k+1} = [m] \backslash i_k$.

The GSSRK method makes use of the Gramian of matrix $A$, that is $G = A A^T$, where each entry can be written as $G_{ij} = \langle A^{(i)}, A^{(j)} \rangle$. It is known that if an equation $A^{(j)} x = b^{(j)}$ is solved by iterate $x^{(k)}$ and if $A^{(i_k)}$ is orthogonal to $A^{(j)}$, that is, if $G_{i_k j} = 0$, then the equation is also solved by the next iterate $x^{(k+1)}$. This means that if $j \notin \mathcal{S}_k$ and if $G_{i_k j} = 0$, then equation $A^{(j)} x = b^{(j)}$ is still solvable by the next iterate $x^{(k+1)}$ and the index $j$ should remain unselectable for one more iteration. However, indices that satisfy $G_{i_k j} \neq 0$ should be reintroduced in the selectable set in each iteration. In summary, the last step of the SSRK algorithm can be written as $S_{k+1} = (S_k \cup \{ j: G_{i_k j} \neq 0 \} ) \backslash \{ i_k \}$.

In terms of performance, the authors observed that: NSSRK and RK are almost identical; GRK outperforms GSSRK, NSSRK, and RK; depending on the data sets, GSSRK can either outperform NSSRK and RK or it can be similar to RK. It is important to note that these conclusions were only made in terms of how fast the error decreases as a function of the number of iterations. There is no analysis of the performance of the algorithms in terms of computation time.

\subsection{Randomized Kaczmarz with Averaging Method} \label{rka_par}

The Randomized Kaczmarz method is difficult to parallelize since it uses sequential updates. Furthermore, just as it was mentioned before, RK does not converge to the least-squares solutions when dealing with inconsistent systems. To overcome these obstacles, the Randomized Kaczmarz with Averaging (RKA) method \cite{moorman2021randomized} was introduced by Moorman, Tu, Molitor, and Needell. It is a block-parallel method that, in each iteration, computes multiple updates that are then gathered and averaged. Let $q$ be the number of threads and $\tau_k$ be the set of $q$ rows randomly sampled in each iteration. In that case, the Kaczmarz step can be written as
\begin{equation} \label{eq:rka_eq}
    x^{(k+1)} = x^{(k)} + \frac{1}{q} \sum_{i \in \tau_k} w_i \frac{b_i - \langle A^{(i)}, x^{(k)} \rangle }{\|A^{(i)}\|^2} {A^{(i)}}^T
\end{equation}
where $w_i$ are the row weights that have a similar function to the relaxation parameter $\alpha_i$ in (\ref{eq:alg}). The projections corresponding to each row in set $\tau_k$ should be computed in parallel. The authors of this method have shown that not only RKA has linear convergence such as RK, but that it is also possible to decrease the convergence horizon for inconsistent systems if more than one thread is used. More specifically, they showed that using the same probability distribution as in RK for row sampling, that is, proportional to the row's squared norm, and using uniform weights $w_i = \alpha$, the convergence can be accelerated such that the error in each iteration satisfies
\begin{equation}
\begin{gathered}
\resizebox{0.9\hsize}{!}{%
    $\EX \| x^* - x^{(k)} \|^2 \leq \sigma_{max} \Big ( \Big ( I - \alpha \frac{A^TA}{\|A\|_F^2} \Big )^2 + \frac{\alpha^2}{q} \Big ( I - \frac{A^TA}{\|A\|_F^2} \Big ) \frac{A^TA}{\|A\|_F^2} \Big ) \| x^* - x^{(0)} \|^2 \:$} \\
    + \frac{\alpha^2}{q} \frac{\|r_{LS}\|^2}{\|A\|_F^2}.
\end{gathered}
\end{equation}
Notice that larger values of $q$ lead to faster convergence and a lower convergence horizon.
The authors also give some insight into how $\alpha$ can be chosen to optimize convergence. In the special case of uniform weights and consistent systems, the optimal value for $\alpha$ is given by:
\begin{equation} \label{eq:best_alpha}
    \alpha^* =
    \begin{cases}
        \mathlarger{\frac{q}{1+(q-1)s_{min}}}, & s_{max} - s_{min} \leq \mathlarger{\frac{1}{q-1}}\\
        \mathlarger{\frac{2q}{1+(q-1)(s_{min}+s_{max})}}, & s_{max} - s_{min} > \mathlarger{\frac{1}{q-1}}
    \end{cases}
\end{equation}
where $s_{min} = \sigma_{min}^2(A) / \|A\|_F^2$ and $s_{max} = \sigma_{max}^2(A) / \|A\|_F^2$.
Although the authors proved that, if the computation of the $q$ projections can be parallelized, RKA can have a faster convergence rate than RK, they did not implement the algorithm using shared or distributed memory, meaning that no results regarding speedups are presented.

\subsection{Summary} \label{table}

\begin{table}[t]
\centering
\caption{Summary of the variations of the Kaczmarz method. In the case of underdetermined consistent systems, $x_{LN}$ refers to the least Euclidean norm solution. For overdetermined systems, $x^*$ refers to the unique solution of consistent systems, and $x_{LS}$ refers to the least-squares solution of inconsistent systems.
}
\begin{tabular}{|l|c|c|c|c|c|}
    \hline
    \hfil \multirow{2}{*}{Name}\hfill & \multirow{2}{*}{Sect.} & \multirow{2}{*}{Year} & \multicolumn{3}{c|}{Convergence} \\ \hhline{~|~|~|-|-|-}
    & & & $x_{LN}$ & $x^*$ & $x_{LS}$ \\ \hline
    Cyclic Kaczmarz (CK) & \ref{sec:CK} & 1937 & $\surd$ & $\surd$ & $\times$ \\ \hline
    Randomized Kaczmarz (RK) & \ref{sec:RK} & 2007 & $\surd$ & $\surd$ & $\times$ \\ \hline
    Randomized Block Kaczmarz (RBK) & \ref{sec:RBK} & 2014 & $\surd$ & $\surd$ & $\times$ \\ \hline
    Randomized Gauss-Seidel (RGS) & \ref{sec:RGS} & 2009 & $\times$ & $\surd$ & $\surd$ \\ \hline
    Randomized Extended Kaczmarz (REK) & \ref{sec:REK} & 2013 & $\surd$ & $\surd$ & $\surd$ \\ \hline
    Randomized Double Block Kaczmarz (RDBK) & \ref{sec:RDBK} & 2015 & $\surd$ & $\surd$ & $\surd$ \\ \hline
    Greedy Randomized Kaczmarz (GRK) & \ref{sec:grk} & 2018 & $\surd$ & $\surd$ & $\times$ \\ \hline
    Non-Repetitive Selectable Set Randomized Kaczmarz (NSSRK) & \ref{ssrk} & 2021 & $\surd$ & $\surd$ & $\times$ \\ \hline
    Gramian Selectable Set Randomized Kaczmarz (GSSRK) & \ref{ssrk} & 2021 & $\surd$ & $\surd$ & $\times$ \\ \hline
    Randomized Kaczmarz with averaging (RKA) & \ref{rka_par} & 2021 & $\surd$ & $\surd$ & $\times$ \\ \hline
\end{tabular}
\label{tab:tab_comp}
\end{table}

We summarize the variations of the original Kaczmarz method analyzed so far in Table~\ref{tab:tab_comp}.


\subsection{Parallel Implementations of the Kaczmarz method} \label{par_kacz}

There are two main approaches to parallelizing iterative algorithms, both of which divide the equations into blocks. The first mode, called block-sequential or block-iterative, processes the blocks in a sequential fashion, while computations on each block are done in parallel. In the second mode of operation, the block-parallel, the blocks are distributed among the processors, and the computations in each block are simultaneous; the results obtained from the blocks are then combined to be used in the following iteration. There are several block-parallel implementations of the original and Randomized Kaczmarz method.

One of the parallel implementations of CK is called Component-Averaged Row Projections (CARP) \cite{gordon2005component}. This method was introduced to solve linear systems derived from partial differential equations (PDEs). It assigns blocks of equations to processors that compute the results in parallel. The results are then merged before the next iteration. This method was developed for both shared and distributed memory and exhibits an almost
linear speedup ratio. However, this implementation was developed for linear systems derived from PDEs that are sparse.

There is a shared memory parallel implementation of the RK method called the Asynchronous Parallel Randomized Kaczmarz (AsyRK) algorithm \cite{liu2014asynchronous} developed specifically for sparse matrices. This algorithm was developed using the HOGWILD! technique~\cite{recht2011hogwild}, a parallelization scheme that minimizes the time spent in synchronization. The idea is to have threads compute iterations in parallel and have them update the solution vector at will. This means that processors can override each other's updates. However, since data is sparse, memory overwrites are minimal and, when do happen, the impact on the solution estimate is negligible. The downsides of this method is that there are restrictions on the number of processors so that the speedup is linear and the matrices must be sparse.

\section{Applications of the Kaczmarz Method} \label{sec:applic_kacz}

We start this Section by describing how the Kaczmarz method can be viewed as a special case of other iterative methods. Then, in Section~\ref{sec:ineq_theory}, we introduce some variations of the method to solve systems of inequalities. Finally, in Section~\ref{sec:kacz_tomo}, we introduced some Kaczmarz-based methods for solving linear systems derived from CT scans.

\subsection{Relationship with Other Iterative Methods} \label{sec:sgd_cimmino}

Chen \cite{chen2018kaczmarz} showed that the Kaczmarz method can be seen as a special case of other row action methods.

A first case worth mentioning is that the Randomized Kaczmarz method can be derived from the Stochastic Gradient Descent (SGD) method using certain parameters. The SGD method is a particular case of the more general Gradient Descent (GD) method. The GD method is a first-order iterative algorithm for finding a local minimum of a differentiable function, the objective function, by taking repeated steps in the opposite direction of the gradient, that is, the direction of the steepest descent. Each iteration of the gradient descent method can be written as
\begin{equation}
    x^{(k+1)} = x^{(k)} - \alpha_k \nabla F(x^{(k)})
\end{equation}
where $\alpha_k$ is the step and $F(x)$ is the objective function to minimize. $F(x)$ is a sum of functions associated with each data entry, represented by $f_i(x)$. In the context of machine learning, the GD method is used to minimize an objective function that usually depends on a very large set of data, which makes it very computationally expensive. The difference between SGD and GD is that in SGD the gradient of the objective function is estimated by using a random subset of data, instead of using the entire data set. The size of the subset is called the minibatch size. If each function, $f_i(x)$, has an associated weight, $w(i)$, that influences the probability of that data entry being chosen, we say that
\begin{equation}
    f_i^{(w)}(x) = \frac{1}{w(i)} f_i(x) \: ,
\end{equation}
and the objective function is written as an expectation
\begin{equation}
    F(x) = \EX^{(w)} [f_i^{(w)}(x)] \: .
\end{equation}
For a more detailed explanation see Section 3 of \cite{needell2014stochastic}. If the minibatch size is one, that is, if we only use one data entry per iteration, the gradient of the objective function can be written as $\nabla F(x) \approx \nabla f_i^{w} (x) = \nabla \frac{1}{w(i)} f_i(x)$, and the algorithm can now be written as
\begin{equation} \label{eq:sgd}
    x^{(k+1)} = x^{(k)} - \frac{\alpha_k}{w(i_k)} \nabla f_{i_k}(x_k)
\end{equation}
where $i_k$ is the chosen data entry in iteration $k$.
Chen \cite{chen2018kaczmarz} showed that the Randomized Kaczmarz method is a special case of SGD using one data entry per iteration. If the objective function is $F(x) = \frac{1}{2} \| Ax-b \|^2$, then $f_i(x) = \frac{1}{2}(b_i \: - \langle A^{(i)}, x^{(k)} \rangle)^2$. The associated weight is the probability of choosing row $i$, that is, $\| A^{(i)} \|^2 / \| A \|_F^2$. Substituting this definition in (\ref{eq:sgd}), we have that
\begin{equation}
    x^{(k+1)} = x^{(k)} + \alpha_k \| A \|_F^2 \frac{b_i \: - \langle A^{(i)}, x^{(k)} \rangle}{\| A^{(i)} \|^2} {A^{(i)}}^T \: ,
\end{equation}
and (\ref{eq:alg}) is recovered if the step is set to $\alpha_k = 1/\| A \|_F^2$.

A second case that we would like to mention is the similarity between the Kaczmarz and the Cimmino method.

The original Cimmino method \cite{cimmino1938cacolo, guida2023reflection} is an iterative method to solve consistent or inconsistent squared systems of equations ($m = n$). A single iteration of the method is made up of $n$ steps. In each step we compute, $y^{(i)}$, the reflection of the previous estimate of the solution with respect to the hyperplane defined by the $i$th equation of the system, given by
\begin{equation}
    y^{(i)} = x^{(k)} + 2 \frac{b_i \: - \langle A^{(i)}, x^{(k)} \rangle}{\| A^{(i)} \|^2} {A^{(i)}}^T \: , \quad \text{with} \quad i = 0, 1, 2, ..., n-1 \: .
\end{equation}
An iteration of the original Cimmino method can then be described as the average of all computed reflections, such that, 
\begin{equation}
    x^{(k+1)} = \frac{1}{n} \sum_{i=0}^{n-1} y^{(i)} = x^{(k)} + \frac{2}{n} \sum_{i=0}^{n-1} \frac{b_i \: - \langle A^{(i)}, x^{(k)} \rangle}{\| A^{(i)} \|^2} {A^{(i)}}^T
\end{equation}
There are several adaptations of the Cimmino method for non-squared matrices and using different relaxation parameters. Here, we present the formulation by Gordon~\cite{gordon2018cimmino}. For an $m\times n$ matrix, an iteration of the Cimmino method with a different relaxation parameter $\alpha_i$ for each row can be written as
\begin{equation} \label{eq:cim_eq_2}
    x^{(k+1)} = x^{(k)} + \frac{1}{m} \sum_{i=0}^{m-1} \alpha_i \frac{b_i \: - \langle A^{(i)}, x^{(k)} \rangle}{\| A^{(i)} \|^2} {A^{(i)}}^T \: .
\end{equation}
Notice that there is a similarity between this adaptation of the Cimmino method and the RKA method (Section~\ref{rka_par}). Expression (\ref{eq:cim_eq_2}) is just the RKA method (Expression~\ref{eq:rka_eq}) if $q = m$, $w_i = \alpha_i$, and if the rows are used ciclicaly.

In practice, it is more common to write the Cimmino method in matrix form such that each iteration is given by:
\begin{equation} \label{eq:cim_matrix}
    x^{(k+1)} = x^{(k)} + A^T D (b - A x^{(k)}) \: ,
\end{equation}
such that $D$ is a diagonal matrix that can be written as $D = \text{diag}(\lambda_i/(m\|A^{(i)}\|^2))$.

In summary, the Cimmino method uses simultaneous reflections of the previous iteration while the Kaczmarz method uses sequential orthogonal projections. The Cimmino method, contrary to the Kaczmarz method, can be easily parallelized. However, the Kaczmarz method converges faster.

\subsection{Adaptations for Linear System of Inequalities} \label{sec:ineq_theory}

The methods described so far solve the linear system of equations presented in (\ref{eq:system}). Another important problem in science and engineering is solving a linear system of inequalities, that is, finding a solution $x$ that satisfies:
\begin{equation} \label{eq:ineq_system}
    Ax \leq b \: .
\end{equation}
These systems can also be classified as consistent or inconsistent if there is or there is not an $x$ that satisfies all the inequalities in (\ref{eq:ineq_system}).

Chen \cite{chen2018kaczmarz} shows how the Kaczmarz method can be a particular case of other iterative methods, some of them for solving linear systems of inequalities. Here we describe two of them.

The relaxation method \cite{agmon1954relaxation, motzkin1954relaxation} introduced by Agmon, Motzkin, and Schoenberg in 1954 is an iterative method that solves the linear system of inequalities (\ref{eq:ineq_system}). An iteration of this method can be computed using
\begin{equation} \label{eq:relax_alg}
    x^{(k+1)} = x^{(k)} + c^{(k)} {A^{(i)}}^T \: , \quad \text{with} \quad i = k \text{ mod } m \: ,
\end{equation}
where $c^{(k)}$ is given by
\begin{equation} \label{eq:par_relax}
    c^{(k)} = min \{ 0, \alpha_i \: \frac{b_i \: - \langle A^{(i)}, x^{(k)}  \rangle}{\|A^{(i)}\|^2} \}.
\end{equation}
Similarly to the Kaczmarz method (see Section~\ref{sec:CK}), the relaxation method has a clear geometric interpretation (if we consider $\alpha_i = 1$): in a given iteration we choose an index $i$ that corresponds to the half-space described by $H_i = \{x : \langle A^{(i)}, x \rangle \leq b_i \}$. More specifically, if the previous iteration $x^{k}$ already satisfies the constraint, the solution estimate is not changed (this corresponds to a negative scale factor, meaning that $c^{(k)} = 0$). Otherwise, the next iteration is given by the orthogonal projection of the previous iteration onto the hyperplane $H_i = \{x : \langle A^{(i)}, x \rangle = b_i \}$. In the original algorithm, rows are chosen cyclically but there are variations with random row selection criteria that admit linear convergence \cite{leventhal2010randomized}.

The relaxation algorithm is guaranteed to converge to a solution that satisfies (\ref{eq:ineq_system}) but, depending on the initial guess of the solution, $x^{(0)}$, the method can converge to different solutions. Hildreth’s algorithm \cite{hildreth1957quadratic} also solves a system of linear inequalities with the addition of converging to the solution closest to the initial guess, $x^{(0)}$. Just like the Kaczmarz method, the random version of Hildreth’s algorithm \cite{jamil2015hildreth} also presents linear convergence to the solution.

\subsection{Application in Image Reconstruction in Computed Tomographies} \label{sec:kacz_tomo}

One of the applications of the Kaczmarz method in the real world is in the reconstruction of images of scanned bodies during a CT scan. 

\subsubsection{From projection data to linear systems}

\begin{figure}[t]
    \centering
    \resizebox{0.6\textwidth}{!}{
    \begin{tikzpicture}
        \draw [draw=red, fill=red, opacity=0.5] plot [smooth cycle] coordinates {(-1,2.5) (-2,1.5) (-3.5,-0.5) (-2,-3) (0,-3) (1,-2) (2,-0.5) (3,0.5) (3,2.5) (1,2.5)};
        \draw[black, -{Latex[length=2mm,width=3mm]}](-16,0)--(16, 0);
        \draw[black, -{Latex[length=2mm,width=3mm]}](0,-16)--(0, 16);
        \node[scale=3, black] at (15, -1) {$x$};
        \node[scale=3, black] at (-1, 15) {$y$};
        \draw[black, -{Latex[length=2mm,width=3mm]}] (-11,7.3)--(1,13.3);
        \draw[black, -{Latex[length=2mm,width=3mm]}](5,-10)--(-5, 10);
        \draw[black, -{Latex[length=2mm,width=3mm]}](6,-9.5)--(-4, 10.5);
        \draw[black, -{Latex[length=2mm,width=3mm]}](7,-9)--(-3, 11);
        \draw[black, -{Latex[length=2mm,width=3mm]}](8,-8.5)--(-2, 11.5);
        \draw[black, -{Latex[length=2mm,width=3mm]}](4,-10.5)--(-6, 9.5);
        \draw[black, -{Latex[length=2mm,width=3mm]}](3,-11)--(-7, 9);
        \draw[black, -{Latex[length=2mm,width=3mm]}](2,-11.5)--(-8, 8.5);
        \coordinate (A) at (-2, 11.5);
        \coordinate (B) at (-3, 11);
        \coordinate (C) at (-4, 10.5);
        \coordinate (D) at (-5, 10);
        \coordinate (E) at (-6, 9.5);
        \coordinate (F) at (-7, 9);
        \coordinate (G) at (-8, 8.5);
        \coordinate (A_new) at ($(A)+(-0.5*0.45,0.5*0.89)$);
        \coordinate (B_new) at ($(B)+(-1*0.45,1*0.89)$);
        \coordinate (C_new) at ($(C)+(-1.5*0.45,1.5*0.89)$);
        \coordinate (D_new) at ($(D)+(-2.5*0.45,2.5*0.89)$);
        \coordinate (E_new) at ($(E)+(-2*0.45,2*0.89)$);
        \coordinate (F_new) at ($(F)+(-1.5*0.45,1.5*0.89)$);
        \coordinate (G_new) at ($(G)+(-0.5*0.45,0.5*0.89)$);
        \draw [draw=blue, very thick] plot [smooth] coordinates {($(-2, 11.5)+(0.5*0.45,0.5*0.89)$) (A_new) (B_new) (C_new) (D_new) (E_new) (F_new) (G_new) ($(-8.5, 8.25)+(0.5*0.45,0.5*0.89)$)};
        \draw[black, dashed] (0,0)--(8,4);
        \draw[black, -{Latex[length=2mm,width=3mm]}] (0:4) arc (0:25:4);
        \node[scale=3, black] at (5,1) {$\theta$};
    \end{tikzpicture}}
    \caption{Spatial domain of a CT scan. The body that is being scanned is represented in red and the projections of each ray are represented in blue. $\theta$ is the rotation angle of the source.}
    \label{fig:ct_im}
\end{figure}
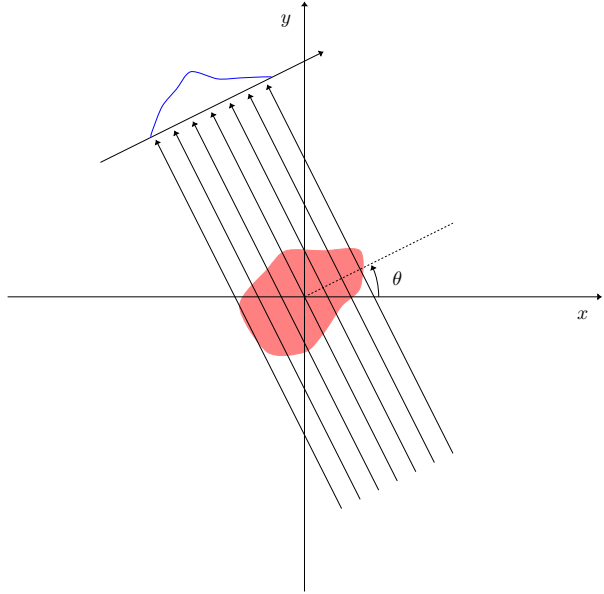

In a CT scan, an X-ray tube (radiation source) and a panel of detectors (target) are spun around the area that is being scanned. For each position of the source and target, defined by an angle $\theta$, the radiation emitted by the X-ray tube is measured by the detectors. After collecting measurements from several positions, the data must be combined and processed to create a cross-sectional X-ray image of the scanned object. We will now explain how the problem of reconstructing the image of a CT scan can be reduced to solving a linear system. Figure~\ref{fig:ct_im} represents the spatial domain of a CT scan.

Let $\mu(r)$ be the attenuation coefficient of the object along a given ray. The damping of that ray through the object is the line integral of $\mu(r)$ along the path of the ray, that is:
\begin{equation} \label{eq:lambert}
    b = \int_{ray} \mu(r) dr \: .
\end{equation}
This is the projection data that is obtained by each detector in the panel. In a CT problem, the goal is to use the projection data, $b$, to find $\mu(r)$, since a reconstructed image of the scanned body can be obtained using the attenuation coefficient of that object. We make the following assumptions to obtain $\mu(r)$ from $b$. First, we assume that the reconstructed image has $L$ pixels on each side, meaning that it is a square of $L^2$ pixels. Secondly, we assume that $\mu(r)$ is constant in each pixel of the image. We can then denote the solution of the CT problem as a vector $x$ with length $L^2$, that contains the attenuation coefficient of the scanned body in every pixel of the image. With these assumptions, we approximate the integral in (\ref{eq:lambert}) using the sum as the simple numerical quadrature, such that
\begin{equation} \label{eq:ray_int}
    b = \sum_{l=1}^{L^2} a_{l} x_l \: ,
\end{equation}
where $a_{l}$ is the length of the ray in pixel $l$. Note that $a_{l}$ is 0 for pixels that are not intersected by the ray. Expression (\ref{eq:ray_int}) corresponds to an approximation of the radiation intensity measured by a single detector. For a given position of source and target, defined by $\theta$, $d$ projections are obtained, where $d$ is the number of detectors in the panel. We can also consider one ray per detector and think of $d$ as the total number of rays for a given position. Considering that measurements are obtained for $\Theta$ angles of rotation, the experimental data comprises $d \times \Theta$ projections. In summary, there are $d \times \Theta$ equations like (\ref{eq:ray_int}). With this setup, we are finally ready to write the CT problem as a linear system.

The system's matrix, $A^{m \times n}$, has dimensions $m = d \times \Theta$ and $n = L^2$. Each row describes how a given ray intersects the image pixels. Calculating the entries of matrix $A$ is not trivial and it can be a very time-consuming process. A vastly used method for computing $A$ is the Siddon method \cite{siddon1985fast}. Since matrix entries are only non-zero for the pixels that are intersected by a given ray, the system's matrices are very sparse. Regarding the other components of the system, the solution $x$ has the pixel's values, and vector $b$ has the projection data.

\subsubsection{Methods for Image Reconstruction}

In the context of image reconstruction during a CT scan, the cyclical Kaczmarz method is also known as the Algebraic Reconstruction Technique (ART) method, first presented by Gordon, Bender, and Herman in 1970 \cite{gordon1970algebraic}. 
Since then, other variations of ART have been developed for systems with special characteristics (for example, systems where each entry of $A$ is either 0 or 1). Now, ART-type methods include methods that use a single data entry in each iteration of the algorithm. For more examples see reference \cite{herman1976iterative}. One variation worth discussing is the incorporation of constraints. A common type of constraints are box constraints that force each entry of the solution to be in between a minimum and maximum value. For example, if it is known that the pixels of the reconstructed image have values between 0 and 1, we can use the $0 \leq x \leq 1$ constraint to decrease the error of the reconstructed image. There are several approaches on how to introduce constraints into our original problem of solving $Ax=b$. Since constraints can be written as inequalities, we could extend the original system to contain equalities and inequalities and use the relaxation method (see Section~\ref{sec:ineq_theory}) for the rows of the system corresponding to the inequalities. In practice, another approach is used: after an iteration of the ART method, the box constraints are applied to every entry in $x$ such that no illegal values are used in the next iteration. This can also be described as the projection of $x^{(k)}$ onto a closed convex set and this variation of the ART method is known as Projected Algebraic Reconstruction Technique (PART) \cite{andersen2014generalized, wu2020projected}. Note that enforcing constraints in each iteration modifies the original linear problem into a non-linear one. Since the original version of ART is based on the CK method, and since the RK method can have a faster convergence than the CK method, the Projected Randomized Kaczmarz (PRK) method \cite{wu2020projected} was introduced to accelerate the convergence of the PART method. This is simply the RK method with an extra step to enforce constraints. Note that, since all these methods are based on variations of the Kaczmarz method for consistent systems, they can not solve inconsistent systems. The authors of PRK solved this problem by developing the Projected Randomzized Extended Kaczmarz (PREK) method, which combines the PART and REK methods to solve inconsistent linear problems with constraints. It is important to note that, in the context of CT scans, there is no consensus on the nomenclature of iterative methods in the literature since there are distinct algorithms with the same name. For instance, there is a parallel implementation of the ART method called PART \cite{gordon2006parallel} that is not related to the projected ART.

Image reconstructions using ART-type methods exhibit salt-and-pepper noise. To solve this problem, the Simultaneous Iterative Reconstruction Technique (SIRT) method was introduced by Gilbert in 1972 \cite{gilbert1972iterative}. Unlike ART-type methods, SIRT uses all equations of the system simultaneously. Once again, ``SIRT'' may refer to different algorithms in literature. Some use SIRT to describe the umbrella of methods that use all equations of the system at the same time. Others use SIRT for the method described in \cite{gilbert1972iterative}.
Similar to the PART method, the Projected Simultaneous Algebraic Reconstruction Technique (PSIRT) method \cite{elfving2012semiconvergence} is a variation of the SIRT algorithm with an extra step in each iteration where the solution estimate solution is projected onto a space that represents certain constraints. SIRT is known to produce fairly smooth images but has slower convergence when compared to ART-type methods.

The Simultaneous Algebraic Reconstruction Technique (SART) method was developed by Andersen and Kak in 1984 \cite{andersen1984simultaneous} as an improvement to the original ART method. It combines the fast convergence of ART with the smooth reconstruction of SIRT. In the ART method, in one iteration, a single data entry that corresponds to one ray for a given scan direction (angle) is processed. In SART, all rays for a given angle are processed at once, such that entry $j$ of solution $x$ in iteration $k+1$ can be computed as follows:
\begin{equation} \label{eq:sart}
    x_{j}^{(k+1)} = x_{j}^{(k)} + \frac{1}{\sum_i A_j^{(i)}} \sum_i \Bigg\{ \frac{b_i \: - \langle A^{(i)}, x^{(k)} \rangle}{\sum_{l}^{N-1} A_l^{(i)}} A_{j}^{(i)} \Bigg\}\: .
\end{equation}
where the sum with respect to $i$ is over the rays intersecting pixel $j$ for a given scan direction. This formula can be generalized to use other relaxation parameters.

\subsubsection{Semi-convergence of Iterative Methods}

Since detectors have measurement errors, it is expected that the raw data from a CT scan is noisy. Let $\overline{b}$ be the error-free vector of constants and $\overline{x}$ be the respective solution. The consistent system can then be written as $A \overline{x} = \overline{b}$. Now suppose that the radiation values given by the detectors are measured with some error $e$. This means that the projections will be given by $b = \overline{b} + e$. The resulting inconsistent system can then be written as $A x_{LS} \approx b$. Note that our goal is not to find the least-squares solution $x_{LS}$, but to find the solution of the error-free case $\overline{x}$.

\begin{figure}[t]
    \centering
    \includegraphics[width=0.5\linewidth]{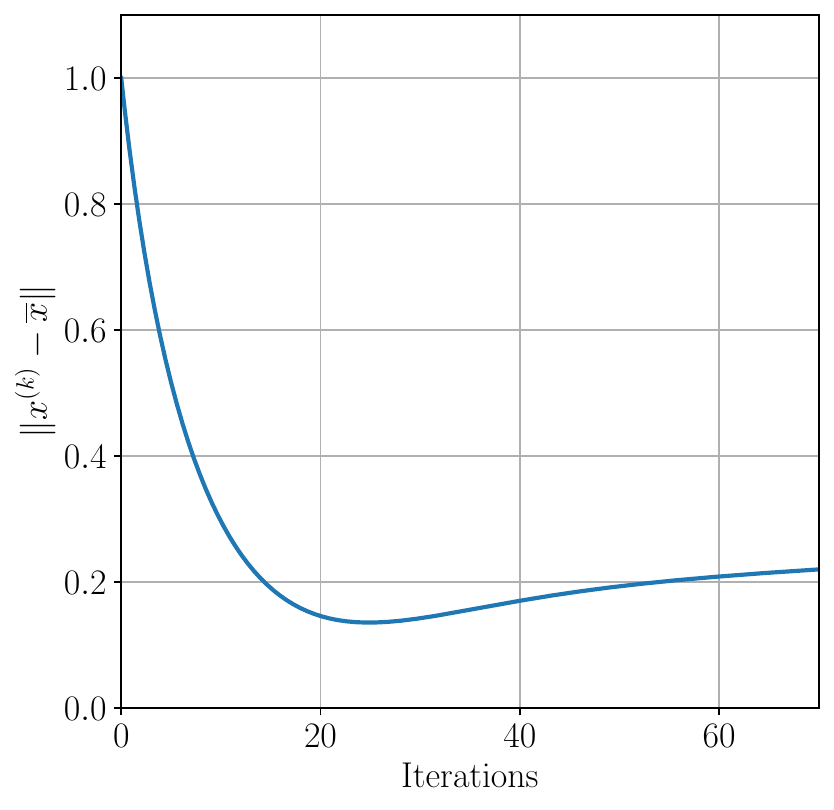}
    \caption{Application on an iterative method to solve problems derived from CT scans. The image shows the evolution of the reconstruction error during several iterations.}
    \label{fig:semi_conv_fig}
\end{figure}

Figure \ref{fig:semi_conv_fig} shows the typical evolution of the reconstruction error, $\|x^{(k)}-\overline{x}\|^2$, for increasing number of iterations when applying an iterative method to solve a linear system from a CT scan. During the initial iterations, the reconstruction error decreases until a minimum is reached. After that, the solution estimate, $x^{(k)}$, will converge to $x_{LS}$ but deviate from the desired value $\overline{x}$. This behaviour is known as semi-convergence~\cite{hansen2018air, hansen2021computed} and the goal is to find the solution estimate at the minimum value of the reconstruction error.


\section{Implementation and Results for the Variations of the Kaczmarz Method}
\label{sec:res_seq}
In this section, we evaluate the performance of implementations of the Kaczmarz method and some of its variants. The section is organized as follows. Sections~\ref{sec:seq_imp} to~\ref{sec:imp_rand} present the implementation details for the considered variants of the Kaczmarz method, as well as a description of the experimental setup. Section~\ref{sec:seq_res} discusses the obtained experimental results.

\subsection{Implemented Methods} \label{sec:seq_imp}

In addition to the original Kaczmarz method (Section~\ref{sec:CK}), the sequential variations from Section~\ref{sec:RK} to Section~\ref{ssrk} were also implemented, with the exception of the block methods in Sections~\ref{sec:RBK} and \ref{sec:RDBK}. These were not implemented since their performance depends strongly on the type of matrix used and they require other external parameters like the partition of the matrix into blocks.

Furthermore, extending on the findings of Needell and Tropp \cite{needell2014paved} that sampling without replacement is an effective row selection criterion (Section~\ref{sec:RBK}), we also implemented the Simple Randomized Kaczmarz Without Replacement (SRKWOR) method. In sampling without replacement in the context of row selection, if a row is selected, it is necessary that all other rows of matrix $A$ are also selected before this row can be selected again. A straightforward technique for sampling rows without replacement is to randomly shuffle the rows of matrix $A$ and then iterate through them in a cyclical manner. We tested two variants of sampling without replacement. In the first, which we called \textit{SRKWOR without shuffling}, the matrix is shuffled only once during preprocessing. In the second method, which we called \textit{SRKWOR with shuffling}, the matrix is shuffled after each cyclical pass. The first variant trades off less randomization for a smaller shuffling overhead than the second variant. Regarding the process of shuffling the rows of the matrix, we consider two possible approaches: the first option is to work with a shuffled array of row indices, a technique we call two-fold indexing; and the second option is to create a new matrix with the desired order. In the first option, the extra memory usage is diminutive but the access time to a given row can be more significant since we first have to access the array of indices and then the respective matrix row. In the second option, we access the matrix directly, while also taking advantage of spatial cache locality: the fact that rows are stored sequentially in memory promotes more cache hits than in the first case. Of these four options, we found that the \textit{SRKWOR without shuffling} using two-fold indexing is the fastest option, meaning that, from this point on, references to the SRKWOR method correspond to this version.

To show that the Kaczmarz method can be a competitive algorithm for solving linear systems, it is worth showing that it can outperform the celebrated Conjugate Gradient method. Therefore, other than the Kaczmarz method and its variants, we will also use the CG \footnote{\url{https://eigen.tuxfamily.org/dox/classEigen_1_1ConjugateGradient.html}} and CGLS \footnote{\url{https://eigen.tuxfamily.org/dox/classEigen_1_1LeastSquaresConjugateGradient.html}} methods from the \textsc{Eigen} linear algebra library \footnote{\url{https://eigen.tuxfamily.org/index.php?title=Main_Page}}. For CG~\footnote{\url{https://eigen.tuxfamily.org/dox/classEigen_1_1DiagonalPreconditioner.html}} we used the default preconditioner (diagonal preconditioner) and both the upper and lower triangular matrices (the parameter \textit{UpLo} was set to \textit{Lower|Upper}). For CGLS~\footnote{\url{https://eigen.tuxfamily.org/dox/classEigen_1_1LeastSquareDiagonalPreconditioner.html}} we also used the default preconditioner (least-squares diagonal preconditioner). Since the CG method can only be used for squared systems, it was necessary to transform the system in (\ref{eq:system}) into a squared system by multiplying both sides of the equation by $A^T$, such that the new system is given by $A^T A x = A^T b$. The computation time of $A^T A$ and $A^T b$ was added to the total execution time of CG.

\subsection{Stopping Criterion} \label{stop_crit_sec}

Since the Kaczmarz method and its variants are iterative methods, it is required to define a stopping criterion. To ensure a fair comparison between the implemented methods and the CG and CGLS methods from the \textsc{Eigen} library, we used the following procedure: first, we determine the number of iterations, $k$, that the different methods take to achieve a given error, that is, $\|x^{(k)}-x^{*}\|^2 < \varepsilon$; then we use those numbers as the maximum number of iterations and measure the execution time of the methods. With this procedure not only do we guarantee that the solutions given by the different methods have similar errors but we also only measure the time spent computing iterations, without taking into account the the computational cost of evaluating the stopping criteria. The value of $\varepsilon$ is defined by the user since it regulates the desired accuracy of the solution. For our simulations, we used $\varepsilon = 10^{-8}$.

\subsection{Data Sets} \label{data_sets_section}

Our simulations will use overdetermined systems, which are more common in real-world problems than underdetermined systems. The Kaczmarz method and its variants were tested for dense matrices using 3 artificially generated datasets. Two data sets with consistent systems were generated: one with contrasting row norms and one with coherent rows. The goal of the first data set, denoted by \textit{DataSet1}, was to evaluate how different row selection criteria perform against matrices with different row norm distributions. The goal of the second data set, denoted by \textit{DataSet2}, was to compare the randomized versions with the original version of the Kaczmarz method for a system similar to the one represented in Figure~\ref{fig:example_coherent}. The goal of the third data set, denoted by \textit{DataSet3}, was to compare the variation of the Kaczmarz method developed for inconsistent systems.

For \textit{DataSet1}, matrix entries were sampled from normal distributions where the average $\mu$ and standard deviation $\sigma$ were obtained randomly. For every row, $\mu$ is a random number between $-5$ and $5$, and $\sigma$ is a random number between $1$ and $20$. The matrix with the largest dimension was generated and smaller-dimension matrices were obtained by ``cropping'' the largest matrix. This keeps some similarities between matrices of different dimensions for comparison purposes. Since these methods will be used in overdetermined systems, matrices can have  $2000$, $4000$, $20000$, $40000$, $80000$, or $160000$ rows and $50$, $100$, $200$, $500$, $750$, $1000$, $2000$, $4000$, $10000$ or $20000$ columns. To guarantee a unique solution for consistent systems, the solution vector $x$ is sampled from a normal distribution with $\mu$ and $\sigma$ using the same procedure as before, and vector $b$ is calculated as the product of $A$ and $x$.

For \textit{DataSet2}, generating matrices with coherent rows can be achieved by having consecutive rows with few changes between them. Just like in the previous data set, we compute the largest matrix and solution and the smaller systems can be constructed by cropping the largest system. Since the analysis of the methods using these data sets will not be as extensive as the first data set, we created fewer matrix dimensions: matrices have 1000 columns and can have 4000, 20000, 40000, 80000, or 160000 rows. The generation of the largest matrix was completed using the following steps: we start by generating the entries of the first row by sampling from a normal distribution $N(2,20)$ - these parameters ensure that the entries are not too similar; we then define the next row as a copy of the previous one, randomly select five different columns, and change the corresponding entries onto new samples of $N(2,20)$; this procedure is repeated until the entire matrix is computed. In summary, we have a random matrix where every two consecutive rows only differ in 5 elements. We computed the maximum angle between any two consecutive rows for the matrix with dimensions $160000 \times 1000$ for \textit{DataSet1} and \textit{DataSet2} and the results are, respectively, $\theta_{max} \approx 2.07$ and $\theta_{max} \approx 0.71$, which confirms that \textit{DataSet2} has much more similar rows (in angle) than \textit{DataSet1}.

Finally, a data set with inconsistent systems, denoted by \textit{DataSet3}, was generated to test the methods that solve least-squares problems (REK and RGS). This was accomplished by adding an error term to the consistent systems from the first data set. Let $b$ and $b_{LS}$ be the vectors of constants of the consistent and inconsistent systems. The latter was defined such that $b_{LS} = b + N(0,1)$. The least-squares solution, $x_{LS}$, was obtained using the CGLS method. The chosen matrix dimensions for this data set are the same as for \textit{DataSet1}.

\subsection{Effect of Randomization} \label{sec:imp_rand}

Since the row selection criterion has a random component in most variants of the algorithm, the solution vector $x$, the maximum number of iterations, and the running time vary with the chosen seed for the random number generator. To get a robust estimate of the number of iterations and execution times, for each input, the algorithm is run several times with different seeds, and the solution $x$ is calculated as the average of the outputs from those runs. We chose to run the algorithm 10 times since this is enough to achieve a percentual standard deviation of $1\%$ regarding the execution time of each run. Furthermore, execution times represented in the results section correspond to the total time of those 10 runs. In every run, the initial estimate guess for the solution, $x^{(0)}$, was set to 0.

\subsection{Results} \label{sec:seq_res}

Simulations were implemented in the \textsc{C++} programming language; its source code and corresponding documentation are publicly available \footnote{Code available here: \url{https://github.com/inesalfe/Review-Seq-Kaczmarz.git}}. All experiments were carried out on a computer with a 3.1GHz central processing (AMD EPYC 9554P 64-Core Processor) and 256 GB memory.

\subsubsection{Randomized Kaczmarz Algorithm} \label{sec:res_seq_rk}

In this section, the simulations for the RK method  used \textit{DataSet1}.

\begin{figure}[t]
    \centering
    \begin{minipage}{.5\textwidth}\centering
    \subfloat[Number of iterations.]{\includegraphics[width=0.95\columnwidth]{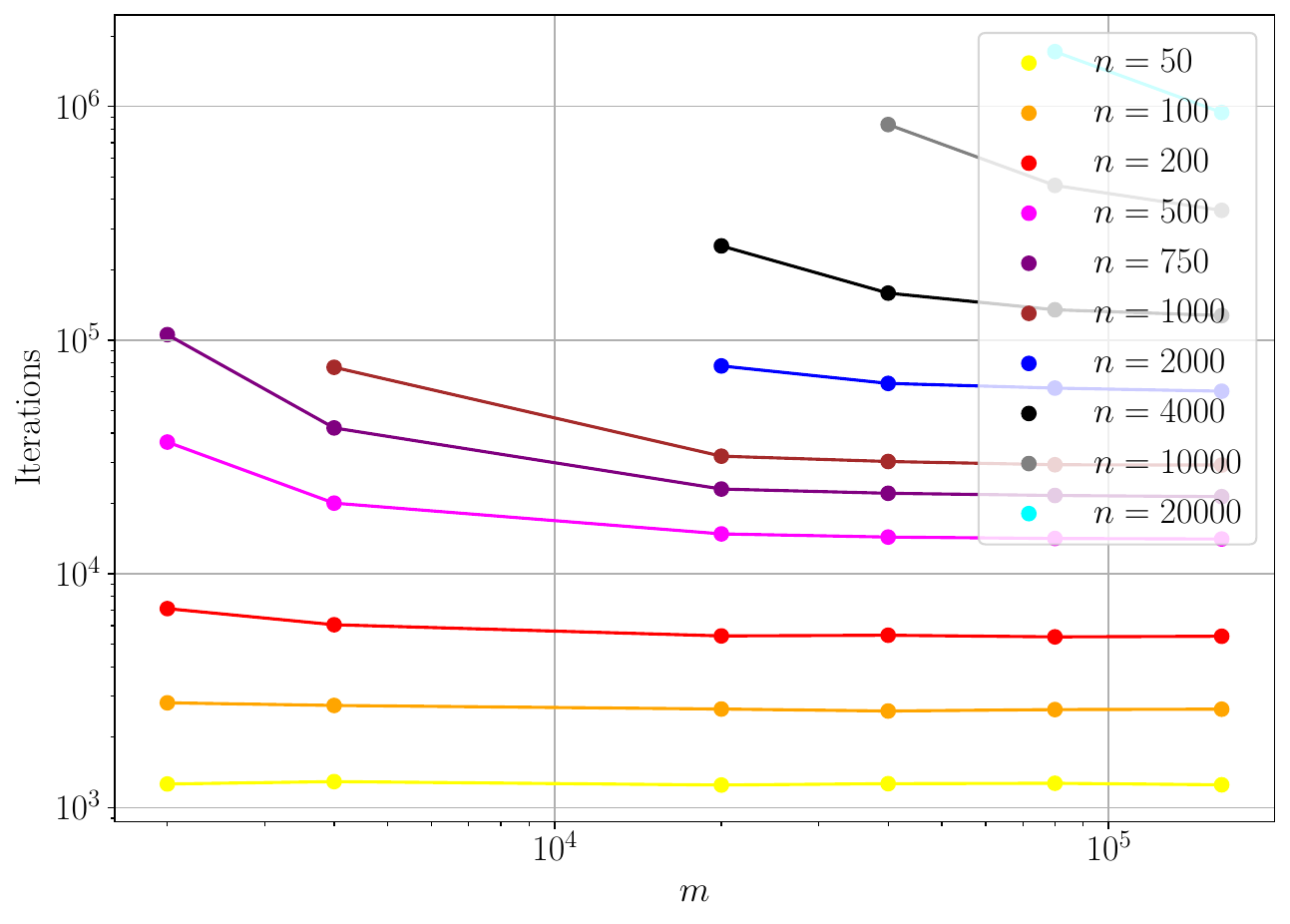}\label{fig:rk_n_1}} \end{minipage}%
    \begin{minipage}{.5\textwidth}\centering
    \subfloat[Execution time.]{\includegraphics[width=0.95\columnwidth]{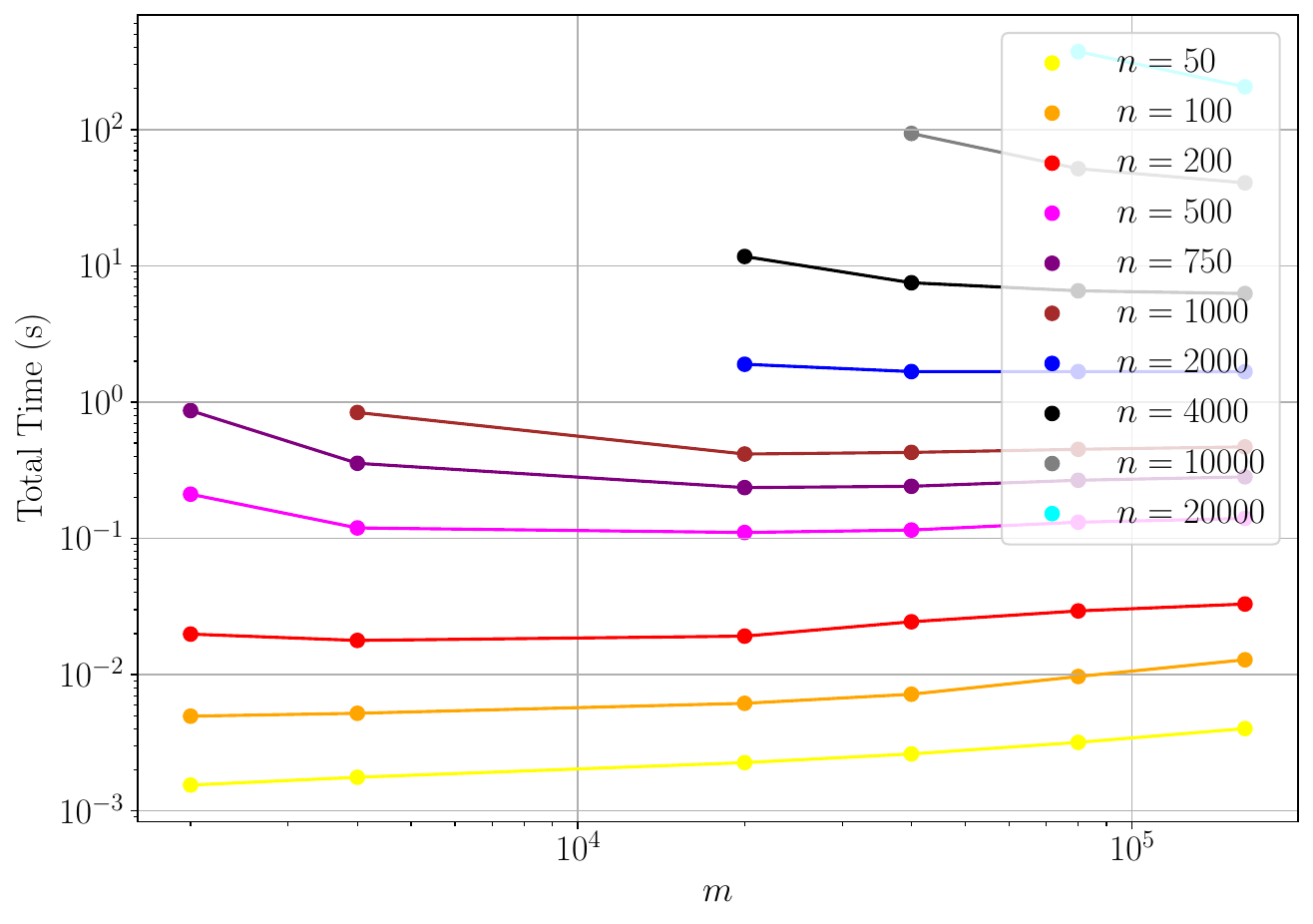}\label{fig:rk_n_2}}
    \end{minipage}%
    \caption{\textit{DataSet1:} Results for the Randomized Kaczmarz algorithm for dense systems using a fixed number of columns and a varying number of rows.}
    \label{fig:rk_n}
\end{figure}

We start with the analysis of the number of iterations and execution time as functions of the number of rows, $m$. From Figure~\ref{fig:rk_n_1} it is clear that, for systems with a fixed $n$, iterations decrease when we increase the number of rows of the system, until they reach a stable value. This is due to the connection between rows and information: for overdetermined systems, more rows translate into more information to solve the system, meaning that increasing $m$ adds restrictions to the system, making it easier to solve and, consequently, lowering the number of iterations. However, there is a point where the number of iterations needed to solve the system is less than the number of rows of the system, after which an increase in $m$ is no longer useful. The initial decrease in iterations is accompanied by a decrease in time, as seen in Figure~\ref{fig:rk_n_2}, after which time starts to increase. This is explained by the fact that, except for the smallest matrix dimension ($2000 \times 50$), all other matrices cannot be stored in cache, and memory access time has a larger impact on runtime when the systems get larger. More specifically, when a row is sampled, the probability that that row is not stored in cache increases for systems with larger $m$, meaning that more time will be spent accessing memory for matrices with a larger number of rows. It is also clear from Figure~\ref{fig:rk_n_1} that the number of iterations increases with $n$. Increasing the number of columns of the system while maintaining $m$ fixed makes the system harder to solve since there are more variables for the same number of restrictions. Figure~\ref{fig:rk_n_2} shows that the total computation time also increases with $n$ due to the correlation with the number of iterations.

\begin{figure}[t]
    \centering    \includegraphics[width=0.5\linewidth]{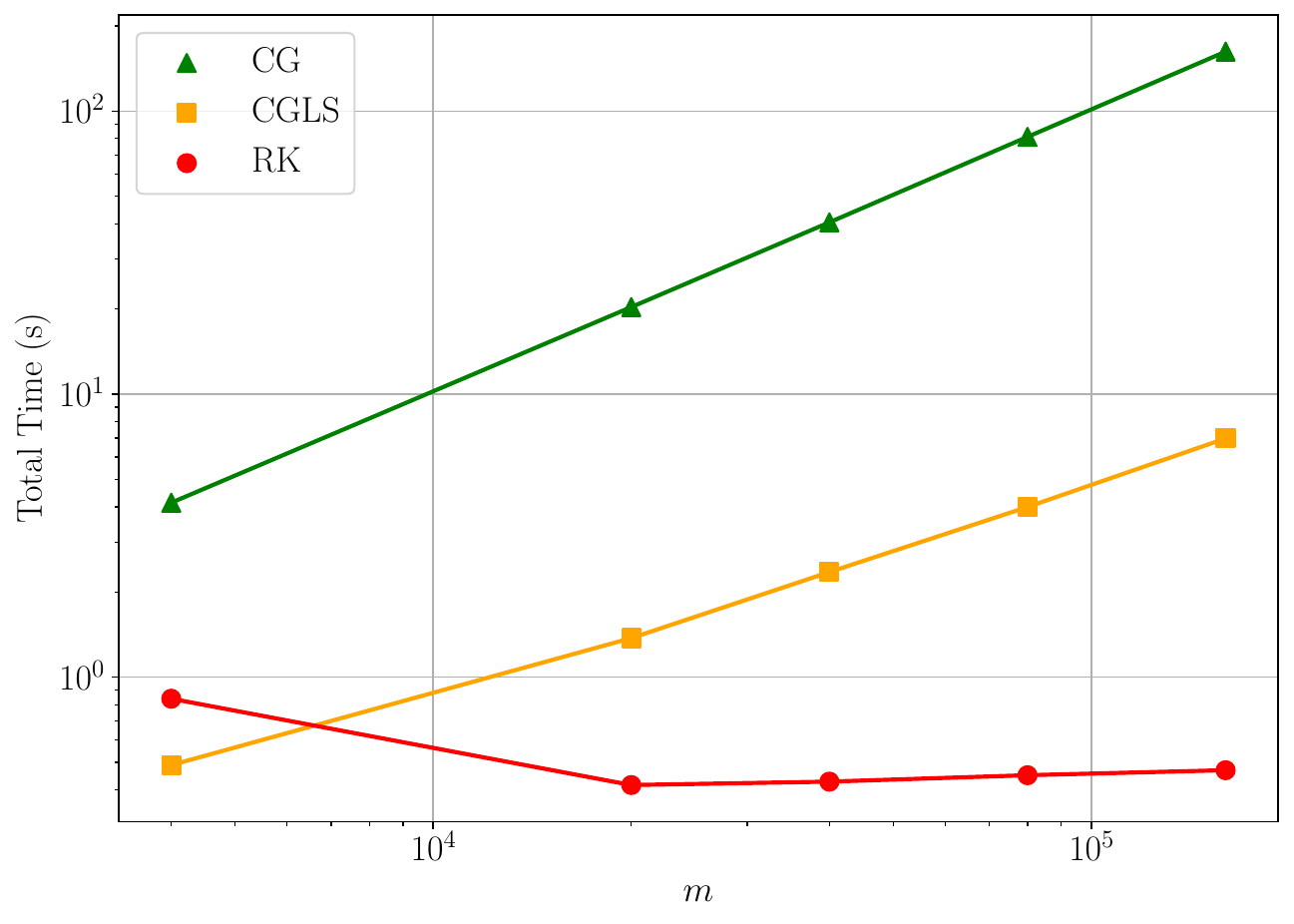}
    \caption{\textit{DataSet1:} Computational time until convergence for RK, CG and CGLS for overdetermined systems with $n=1000$.}
    \label{fig:cg}
\end{figure}

To finish the analysis of RK, we show how this method can outperform CG and CGLS. From Figure~\ref{fig:cg} we can conclude that, regardless of matrix dimension, CGLS is a faster method than CG for solving overdetermined problems. This is to be expected since CGLS is an extension of CG to non-square systems and the execution time of CG takes into account computing $A^T A$. We can see that, except for the smallest system, RK is faster than CGLS. In summary, RK should be used to the detriment of CG and/or CGLS for very overdetermined systems. However, it is important to note that, for the CG method, a very large portion of the execution time is used to transform the system in (\ref{eq:system}) into a squared system. In fact, if we were to analyze the time spent only in iterations of the CG method, this would be faster than RK. Nonetheless, we are working with overdetermined systems and the computation time of $A^T A$ and $A^T b$ cannot be ignored.


\subsubsection{Variants of the Kaczmarz Algorithm for Consistent Systems} \label{sec:var_consist}

\begin{figure}[t]
    \centering
    \begin{minipage}{.5\textwidth}\centering
    \subfloat[Number of iterations. Note that the results for GSSRK, NSSRK, and RK are overlapped.]{\includegraphics[width=0.95\columnwidth]{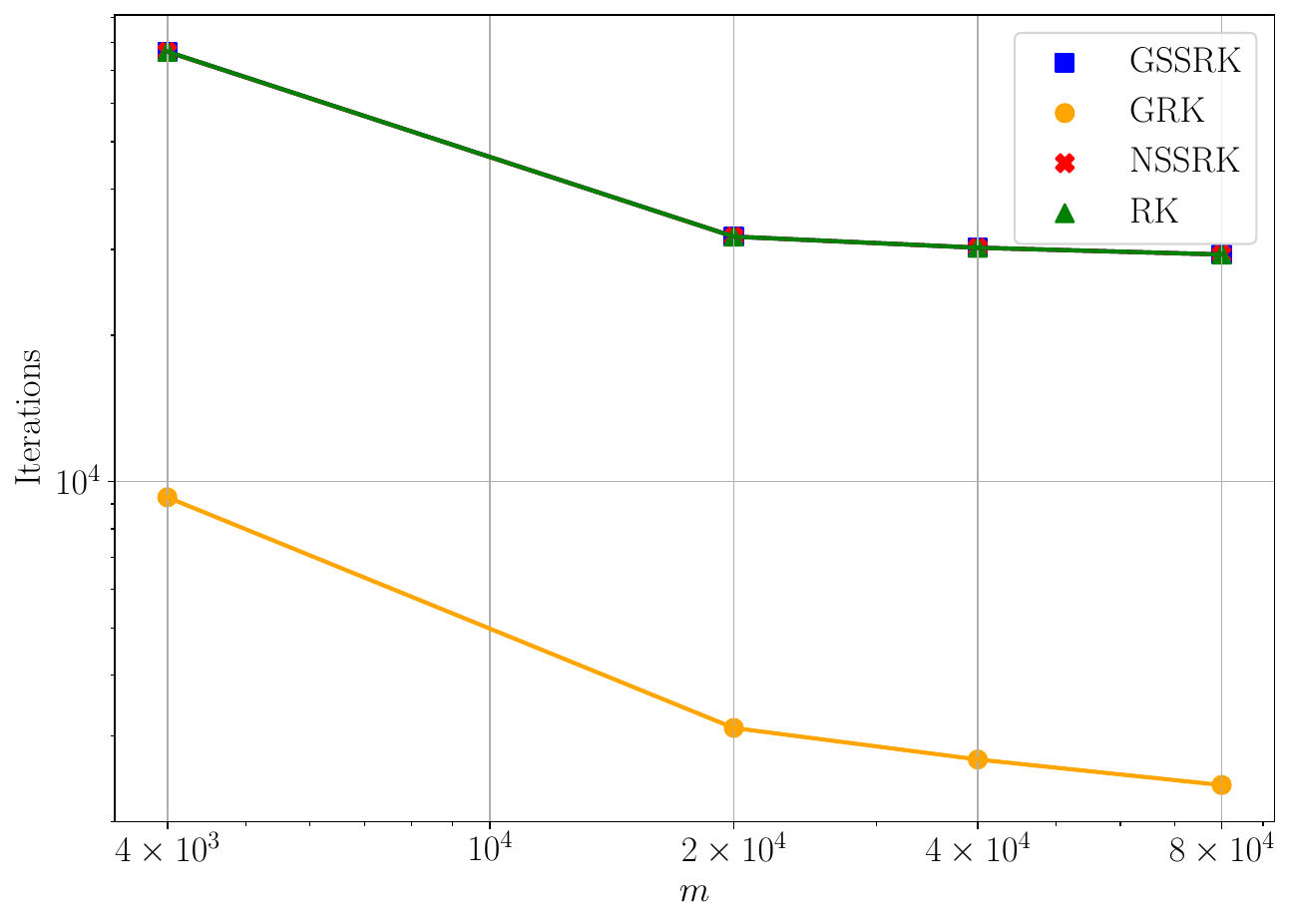}\label{fig:all_ori_n_1}}
    \end{minipage}%
    \begin{minipage}{.5\textwidth}\centering
    \subfloat[Execution time.]{\includegraphics[width=0.95\columnwidth]{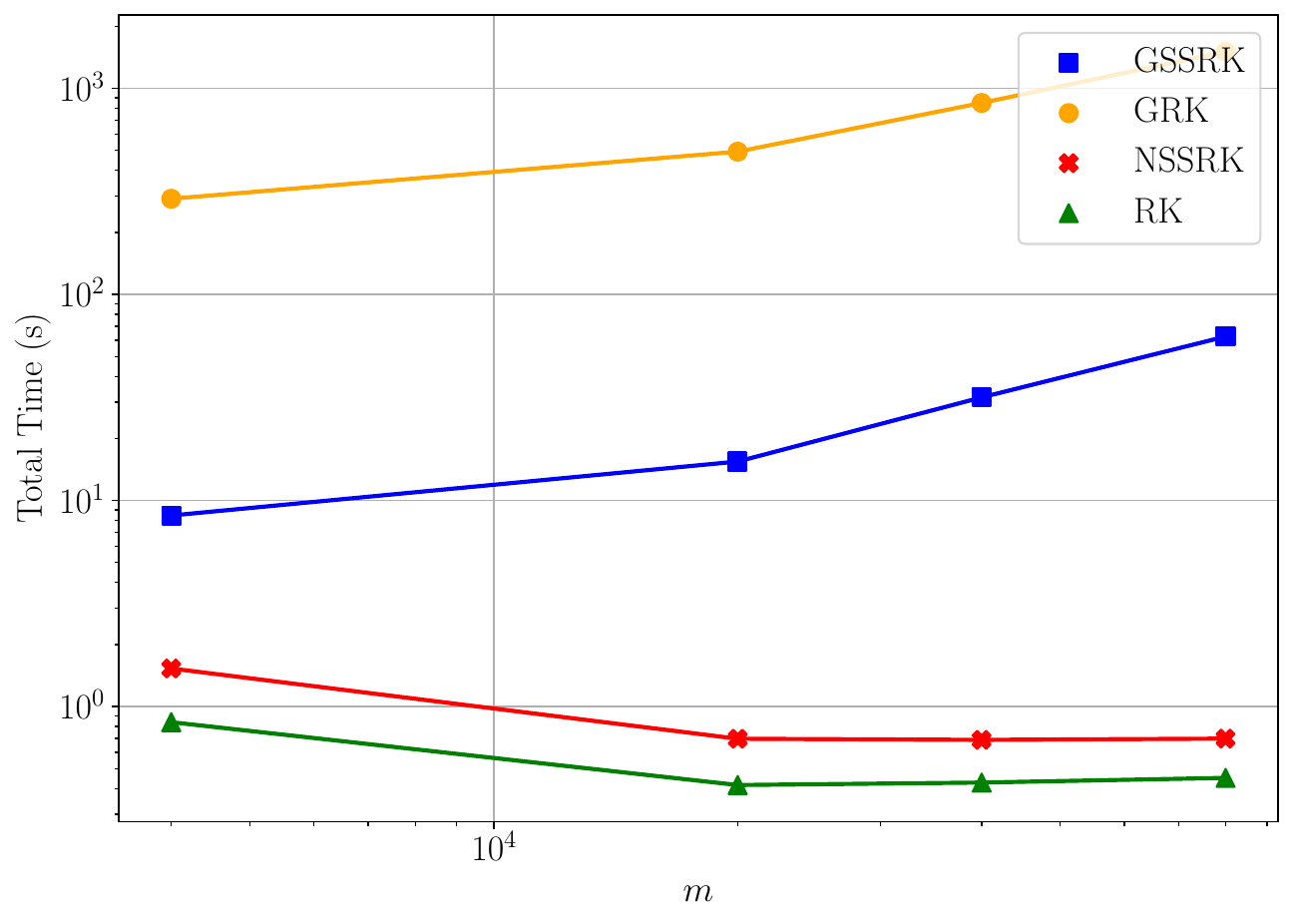}\label{fig:all_ori_n_2}}
    \end{minipage}%
    \caption{\textit{DataSet1:} Results for some variants of the Kaczmarz algorithm for dense systems using a fixed number of columns $n = 1000$ and a varying number of rows.}
    \label{fig:all_ori_n}
\end{figure}

We now compare several variants of the Kaczmarz method between themselves. In this section, we will analyze the results using a fixed number of columns, starting with \textit{DataSet1}. The results for other values of $n$ are similar so we will not discuss them here. We start with the methods introduced in Section~\ref{sec:theory} that, like RK, select rows based on their norms. Figure~\ref{fig:all_ori_n_1} contains the evolution of the number of iterations while the right plot contains the total execution time until convergence. From the number of iterations we can conclude that the GRK method has the most efficient row selection criterion. However, this does not translate into a lower execution time, as observed in Figure~\ref{fig:all_ori_n_1}, meaning that each individual iteration is more computationally expensive than individual iterations of the other methods. This is due to two factors: first, there is the need to calculate the residual in each iteration; second, since rows are chosen using a probability distribution that relies on the residual, and since the residual changes in each iteration, there is the need to update, in each iteration, the discrete probability distribution that is used to sample a single row. For other methods, the discrete probability distribution used for row selection depends only on matrix $A$, and therefore, is the same for all iterations of the algorithm. The RK, GSSRK, and NSSRK methods have an indistinguishable number of iterations. In terms of time (Figure~\ref{fig:all_ori_n_2}), NSSRK and RK have similar performance while GSSRK is a slower method. This happens since, for dense matrices, it is very rare that rows of matrix $A$ are orthogonal and the selectable set usually corresponds to all the rows of the matrix - this means that a lot of time is spent checking for orthogonal rows and very few updates are made to the selectable set. It is also normal for NSSRK to be slightly slower than RK since we are checking if the selected row in a given iteration is the same as the row sampled in the previous one and, for random matrices with many rows, it is unlikely for the same row to be sampled twice in consecutive iterations.

\begin{figure}[b]
    \centering
    \begin{minipage}{.5\textwidth}\centering
    \subfloat[Number of iterations. Note that the results for NSSRK and RK are overlapped.]{\includegraphics[width=0.95\columnwidth]{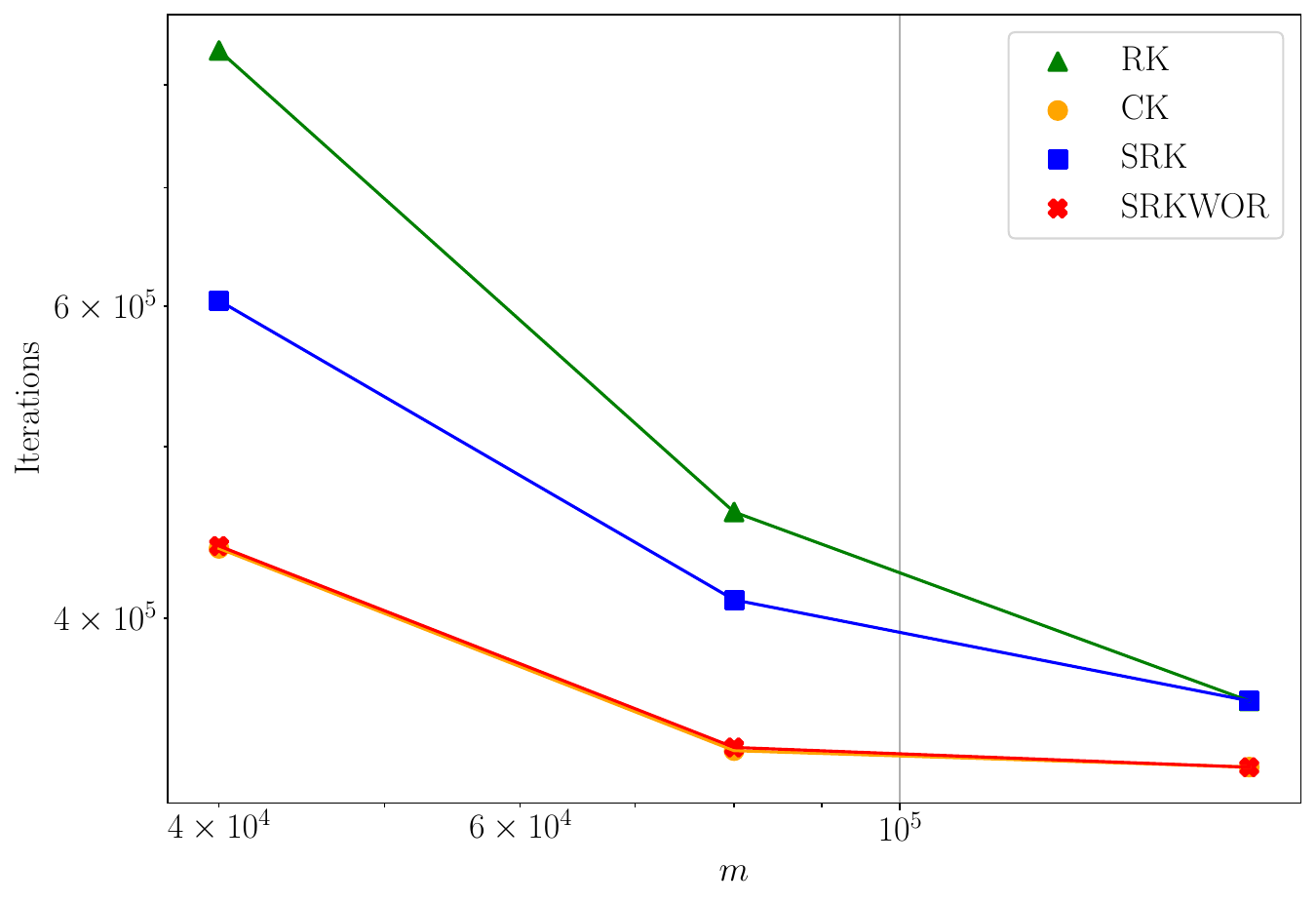}\label{fig:all_exp_n_1}}
    \end{minipage}%
    \begin{minipage}{.5\textwidth}\centering
    \subfloat[Computational time.]{\includegraphics[width=0.95\columnwidth]{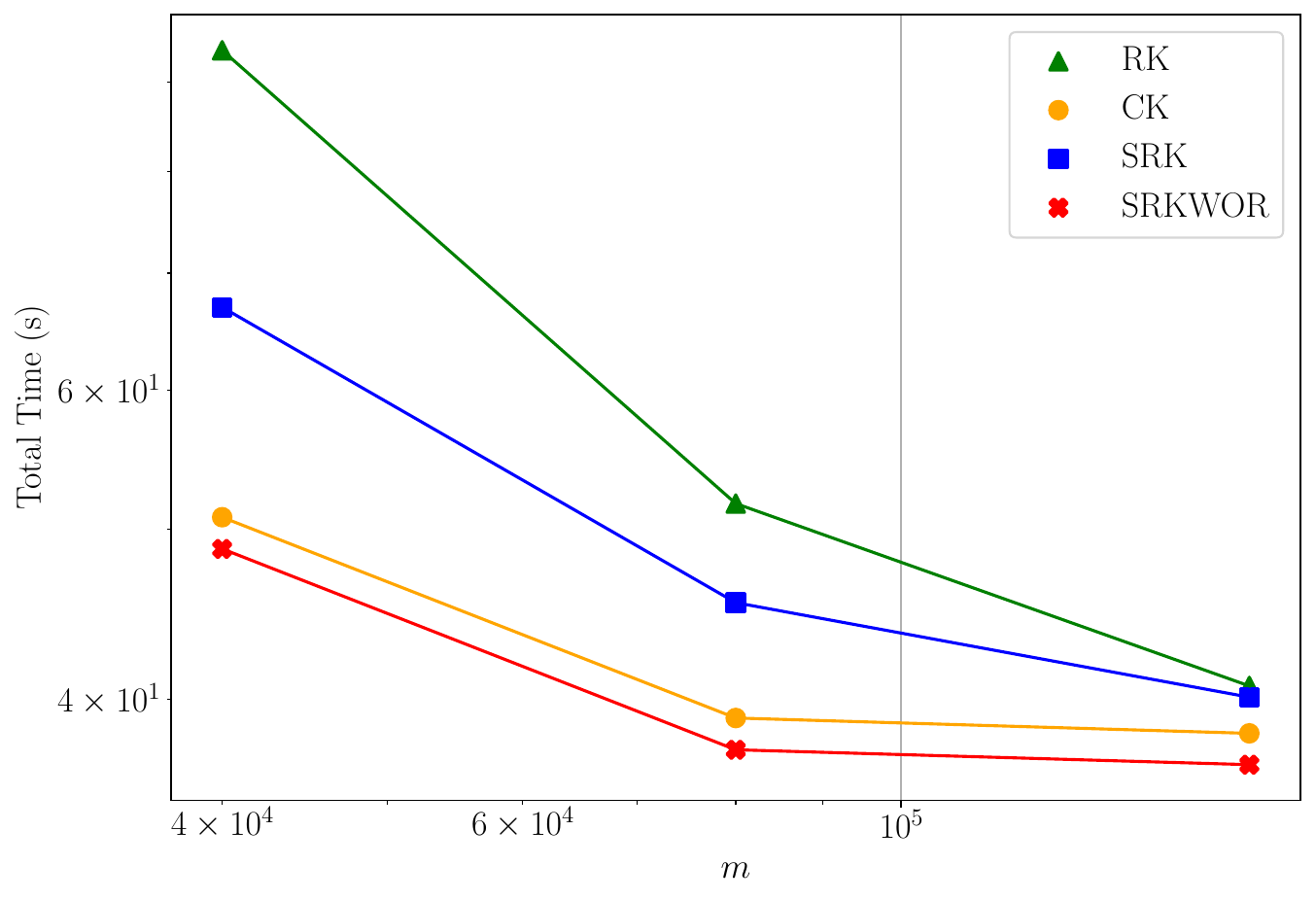}\label{fig:all_exp_n_2}}
    \end{minipage}%
    \caption{\textit{DataSet1:} Results for some variants of the Kaczmarz algorithm for systems with $n = 10000$.}
    \label{fig:all_exp_n}
\end{figure}

When Strohmer and Vershynin in \cite{Strohmer2007ARK} introduced the RK method, they compared its performance to the CK and the SRK methods. Here, we make the same analysis with the addition of SRKWOR. The results for the number of iterations and computational time are presented in Figure~\ref{fig:all_exp_n}. The first observation to be made is that CK yields very similar results to SRKWOR, which is expected when working with random matrices. The SRK method, for smaller values of $m$, exhibits a lower number of iterations and time than the RK method, which shows that, for this data set, sampling rows with probabilities proportional to their norms is not a better row sampling criterion than using a uniform probability distribution. Furthermore, CK and SRKWOR outperform RK and SRK for most dimensions both in terms of iterations and time. These results show that sampling without replacement is indeed an efficient way to choose rows - this may happen since, when using a random approach with replacement, the same row can be chosen many times which makes progress to the solution slower.

\begin{figure}[b]
    \centering
    \begin{minipage}{.5\textwidth}\centering
    \subfloat[Number of iterations.]{\includegraphics[width=0.95\columnwidth]{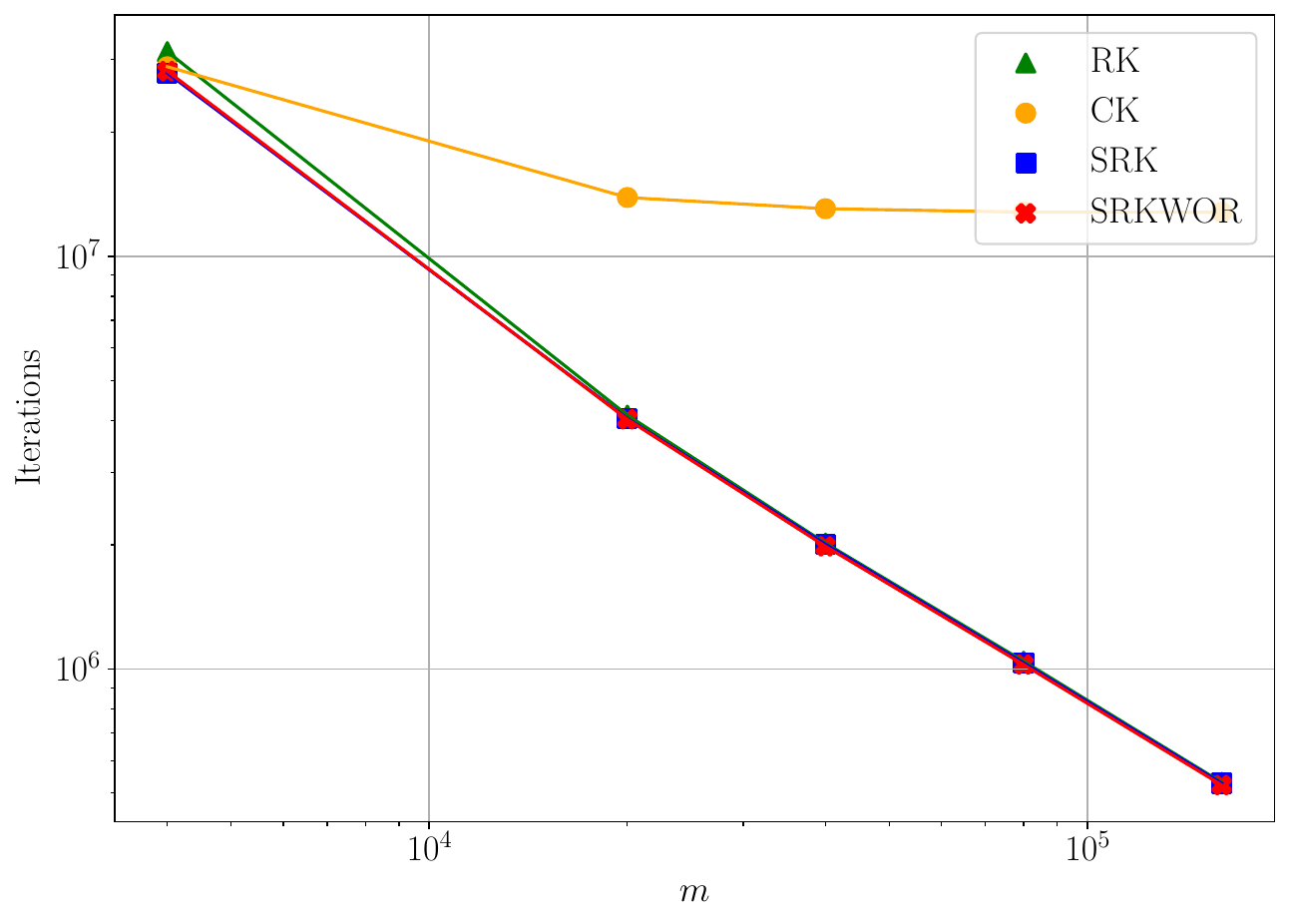}\label{fig:coherent_1}}
    \end{minipage}%
    \begin{minipage}{.5\textwidth}\centering
    \subfloat[Execution time.]{\includegraphics[width=0.95\columnwidth]{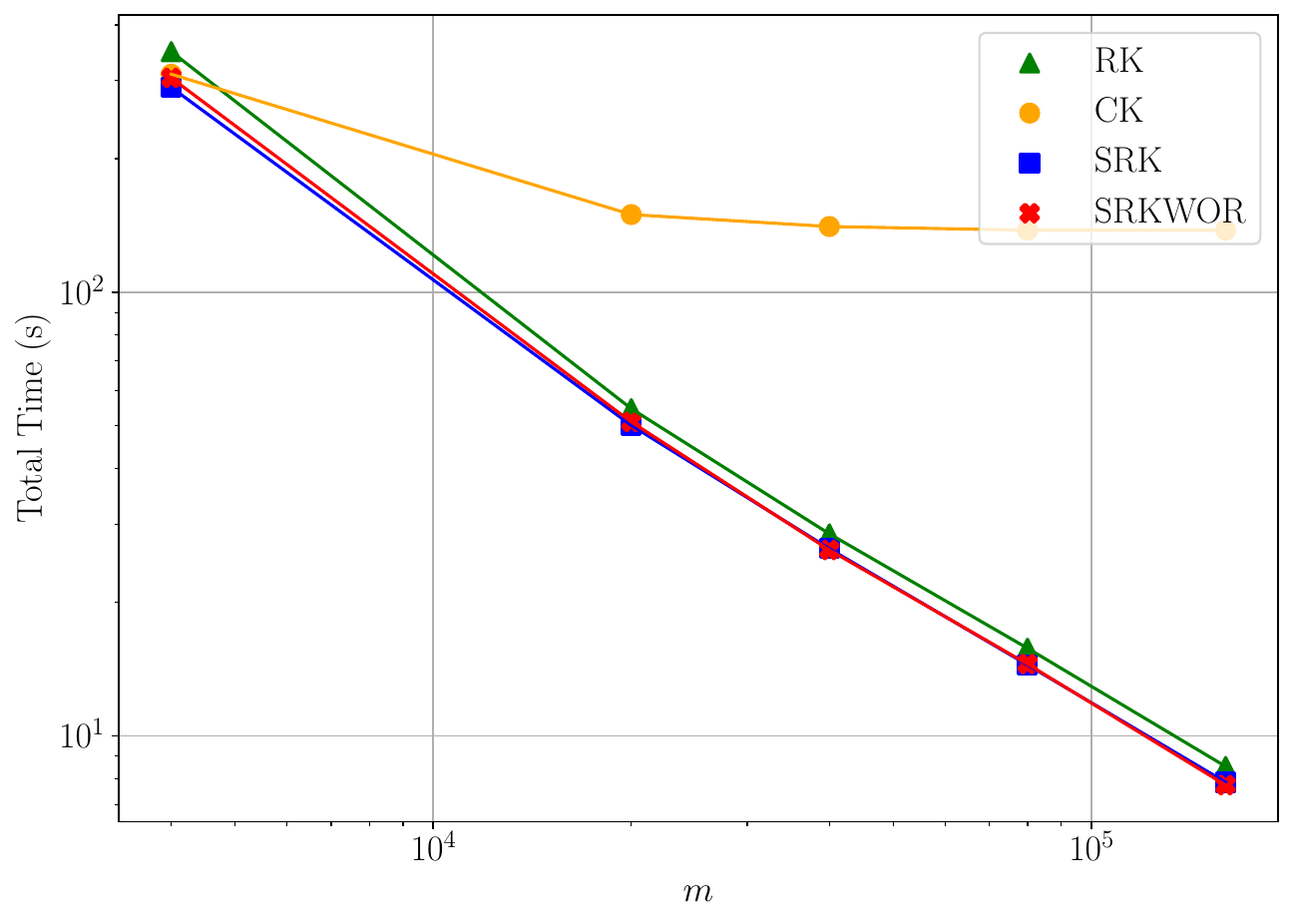}\label{fig:coherent_2}}
    \end{minipage}%
    \caption{\textit{DataSet2:} Results for some variants of the Kaczmarz algorithm for the coherent dense data set using a fixed number of columns $n = 1000$ and a varying number of rows.}
    \label{fig:coherent}
\end{figure}

Contrary to the intuition provided by Strohmer and Vershynin \cite{Strohmer2007ARK}, the RK-inspired method from Figure~\ref{fig:all_exp_n} does not outperform CK. As Wallace and Sekmen~\cite{wallace2014deterministic} refer, choosing rows in a random way should only outperform the CK method for matrices where the angle between consecutive rows is very small, that is, highly coherent matrices. To confirm this, we also compare the results of the methods used in Figure~\ref{fig:all_exp_n} for \textit{DataSet2}. Figure~\ref{fig:coherent} shows that, for highly coherent matrices, CK does have a slower convergence compared to RK and its random variants.

We finish this section with some simulations of the Kaczmarz method using quasirandom numbers, described in Section~\ref{sec:quasirand_intro}. So far, it seems that, for the first dense data set, the fastest method is SRKWOR. For this reason, we compare sampling rows using the quasirandom numbers generated with the Sobol and Halton sequences (SRK-Sobol and SRK-Halton) with the RK and SRKWOR methods.

\begin{figure}[t]
    \centering
    \begin{minipage}{.5\textwidth}\centering
    \subfloat[Number of iterations.]{\includegraphics[width=0.95\linewidth]{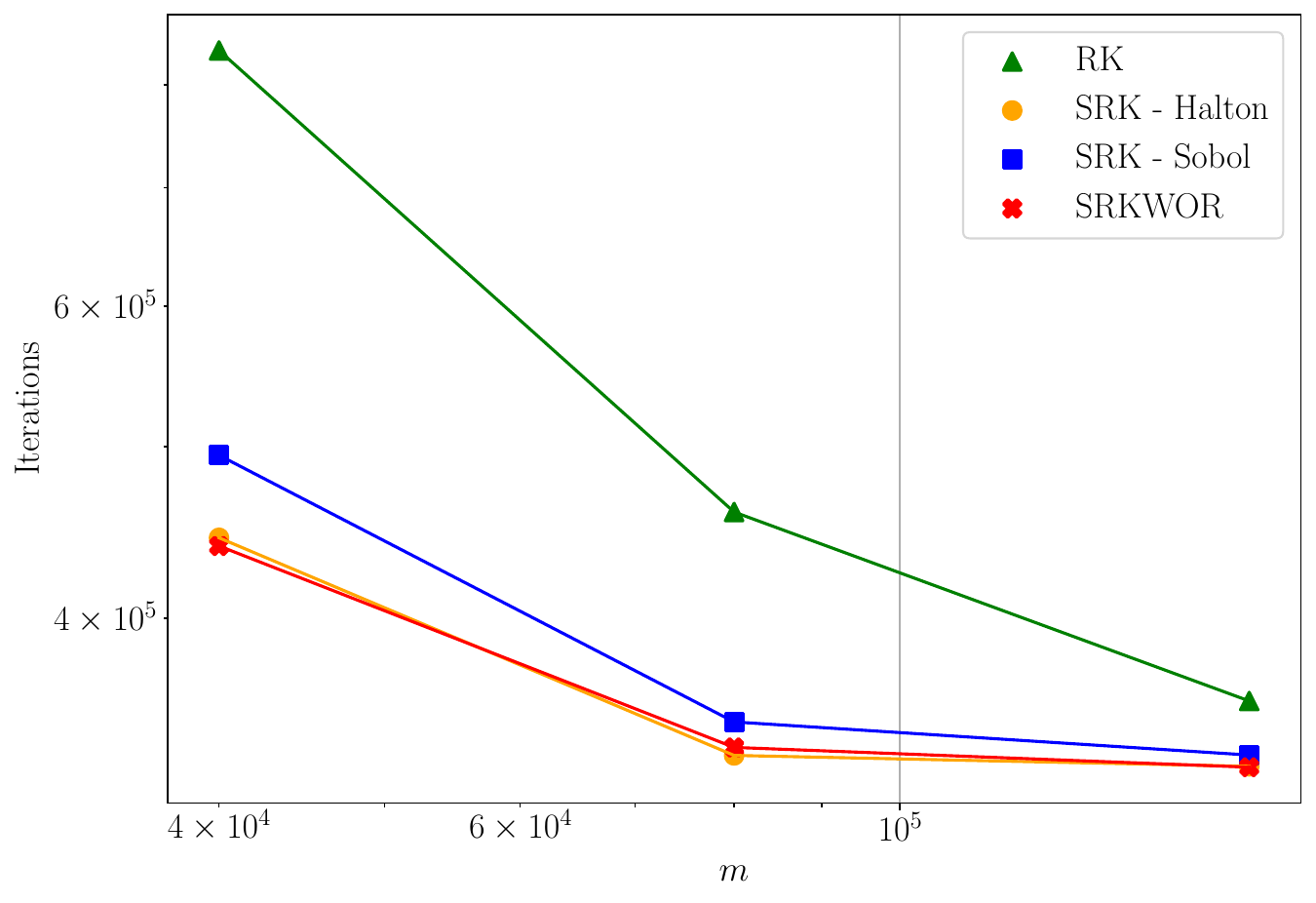}\label{fig:quasirand_n_1}}
    \end{minipage}%
    \begin{minipage}{.5\textwidth}\centering
    \subfloat[Execution time.]{\includegraphics[width=0.95\linewidth]{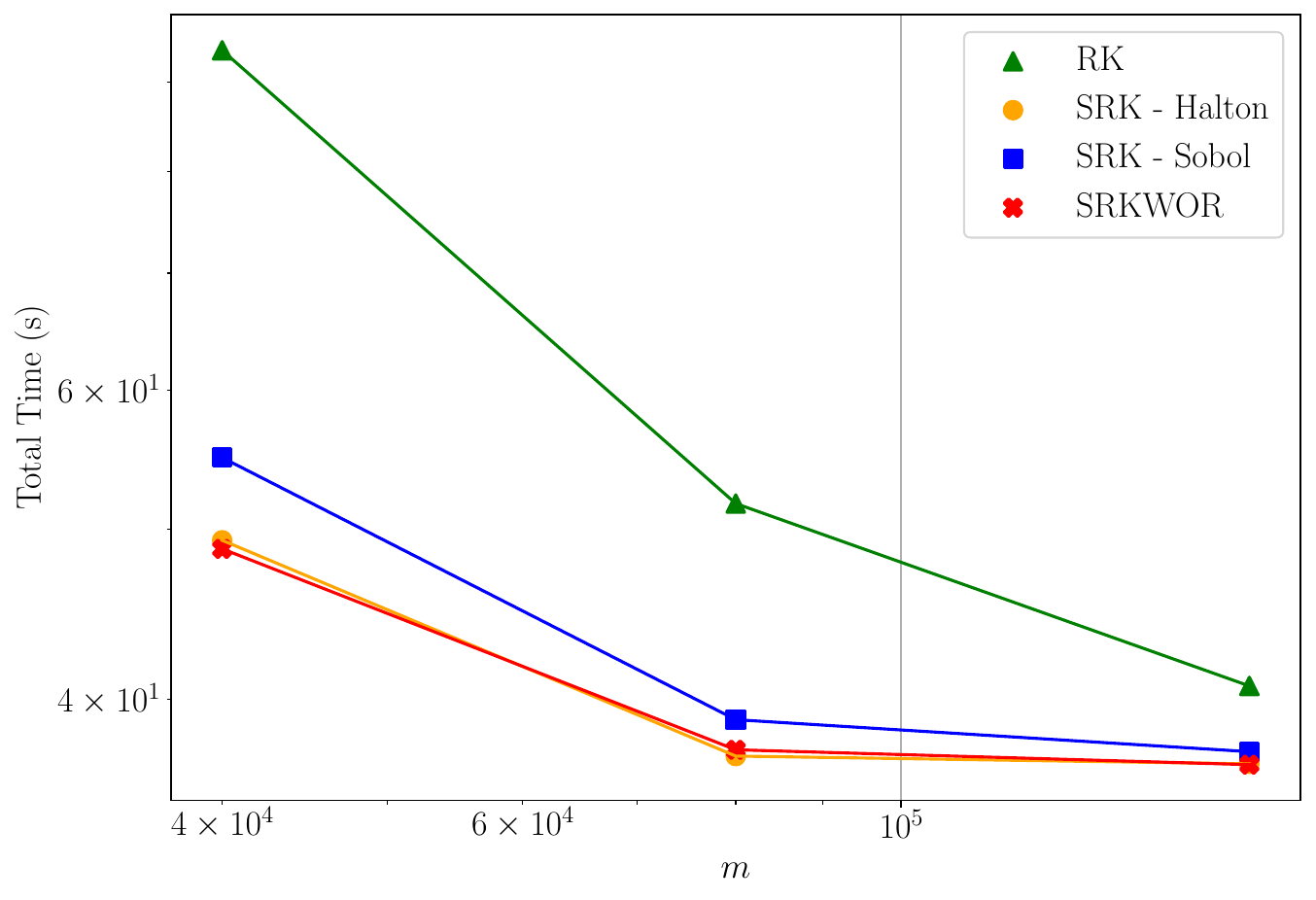}\label{fig:quasirand_n_2}}
    \end{minipage}%
    \caption{\textit{DataSet1:} Results for some variants of the Kaczmarz algorithm for dense systems using a fixed number of columns $n = 10000$ and a varying number of rows.}
    \label{fig:quasirand_n}
\end{figure}

From the results presented in Figure~\ref{fig:quasirand_n}, we can draw two main conclusions: first, quasirandom numbers can outperform RK in both iterations and time; second, the results for the Sobol sequence are very similar to the SRKWOR method, again in terms of iterations and time. In summary, sampling rows randomly without replacement or using quasirandom numbers generates similar results, both faster than RK.

From the analysis of the variants of the Kaczamrz algorithm for consistent systems, we can conclude that, for the data sets that we used, the only methods capable of outperforming the Randomized Kaczmarz method are SRK (using random and quasirandom numbers) and SRKWOR. From all these methods the SRKWOR and SRK-Halton methods appear to be the fastest.

Since many of these variations use row norms to sample rows, we repeated the simulations from this section by changing  \textit{DataSet1} such that rows have similar norms in order to assess if there is a difference in the comparative performance of these methods. However, the results were identical to the ones presented here with one exception: RK and SRK have similar iterations and time since, by having rows with similar norms, the probability distribution used by RK will be similar to a uniform probability distribution, used by SRK.

\begin{figure}[b]
    \centering
    \includegraphics[width=0.5\textwidth]{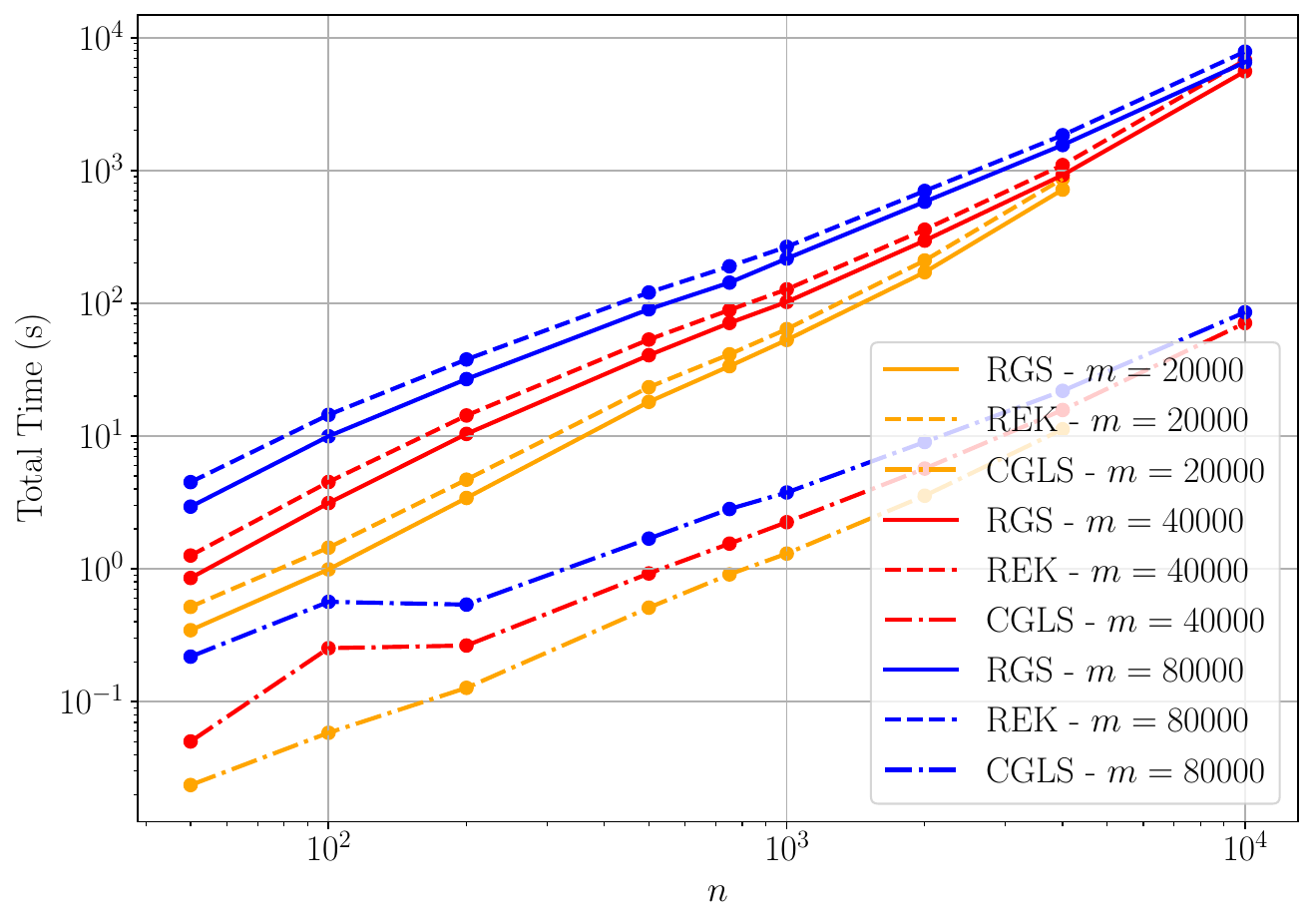}
    \caption{\textit{DataSet3:} Execution time for REK, RGS, and CGLS for dense inconsistent systems using a fixed number of rows and a varying number of columns.}
    \label{fig:ls_inconsist}
\end{figure}

\subsubsection{Variants of the Kaczmarz Algorithm for Least-Squares Systems} \label{sec:var_ls}

In this section, we discuss the results for the REK and RGS methods. In Section~\ref{sec:res_seq_rk}, we compared the RK method with CG and CGLS and concluded that the CGLS is faster than the CG method for overdetermined systems. For this reason, we will also compare REK and RGS with CGLS. The execution time of the three methods, presented in Figure~\ref{fig:ls_inconsist}, shows that, regardless of the dimension of the problem, CGLS is faster than REK and RGS. Between REK and RGS, RGS is slightly faster than REK. In conclusion, REK and RGS are not the most suitable methods for solving least-squares problems and the CGLS method is a faster alternative. Table~\ref{tab:tab_summary_seq} has a summary of the results obtained in this section.

\begin{table}[t]
\centering
\caption{Summary of the results obtained in this section for the CG/CGLS methods and the Kaczmarz method and its sequential variations. Note that we only studied overdetermined systems.}
\begin{tabular}{|c|c|c|c|} \hline
\multicolumn{3}{|c|}{Type of System} & Best Method(s) \\ \hline
\multirow{3}{*}{Consistent} & \multicolumn{2}{c|}{$m \sim n$} & CGLS \\ \hhline{~|-|-|-}
& \multirow{2}{*}{$m >> n$} & Non-coherent & CK / SRKWOR / Quasirandom \\  \hhline{~|~|-|-}
& & Coherent & SRKWOR / Quasirandom \\ \hline
\multicolumn{3}{|c|}{Inconsistent} & CGLS \\ \hline
\end{tabular}
\label{tab:tab_summary_seq}
\end{table}

\section{Conclusion}
\label{sec:conclusion}
In this paper, we surveyed several algorithms based on the Kaczmarz method by describing their row selection criteria and implementing them to evaluate their performance in terms of iteration and execution time for overdetermined dense systems. We started by showing that, for consistent systems, the Randomized Kaczmarz (RK) method can outperform the Conjugate Gradient (CG) and Conjugate Gradient for Least-Squares (CGLS) methods, especially in the case of highly overdetermined systems. For overdetermined systems with $m \sim n$, the CGLS method is a better alternative to RK (and CG). Furthermore, the original Kaczmarz method (CK), the Simple Randomized Kaczmarz Without Replacement (SRKWOR) method, and the Simple Randomized Kaczmarz with quasirandom numbers (SRK-Halton and SRK-Sobol) methods are faster than RK. However, for a special class of consistent systems for which the angle between consecutive rows is small, the randomized versions of the Kaczmarz method are faster than CK. For inconsistent systems, although the Randomized Gauss-Seidel (RGS) method is faster than the Randomized Extended Kaczmarz (REK) method, CGLS is significantly faster than both these methods.

\bmhead{Acknowledgements}

This work was supported by national funds through Fun\-da\-ção para a Ciência e a Tecnologia under projects URA-HPC PTDC/08838/2022 and   UIDB/50021/2020 (DOI:10.54499/UIDB/50021/2020).
JA was funded by  Ministerio de Universidades and specifically the requalification program of the Spanish University System 2021-2023 at the Carlos III University.

\section*{Declarations}

\subsection*{Funding}

This work was supported by national funds through Fun\-da\-ção para a Ciência e a Tecnologia under projects URA-HPC PTDC/08838/2022 and   UIDB/50021/2020 (DOI:10.54499/UIDB/50021/2020).

JA was funded by  Ministerio de Universidades and specifically the requalification program of the Spanish University System 2021-2023 at the Carlos III University.
 
\subsection*{Competing interests}

The authors have no competing interests to declare that are relevant to the content of this article.

\subsection*{Ethics approval}

The authors declare that there are no relevant ethical issues with respect to the content of this article.

\subsection*{Consent for publication}

The authors declare that all consent required for the publication of this paper has been granted.

\subsection*{Data / Code availability}

Data was generated using the publicly available code in \url{https://github.com/inesalfe/Review-Seq-Kaczmarz.git}.

\subsection*{Author contribution}

All authors contributed equally to this work.

\bibliography{references}

\end{document}